\documentclass[12pt,reqno]{amsart}
\usepackage{geometry}
\usepackage[latin1]{inputenc}
\usepackage[italian,english]{babel}
\usepackage{amsmath, amsfonts, amsthm}
\usepackage{graphicx}
\usepackage{pgfplots}
\usepackage{paralist}
\usepackage{hyperref,comment}
\numberwithin{equation}{section}
\newtheorem{thm}{Theorem}[section]

\newtheorem{prop}[thm]{Proposition}
\newtheorem{exa}[thm]{Example}
\newtheorem{defn}[thm]{Definition}
\newtheorem{rem}[thm]{Remark}

\newcommand{\R}{\mathbb{R}}
\newcommand{\N}{\mathbb{N}}

\newcommand{\argmin}{\arg\min}
\newcommand{\eps}{\varepsilon}

\DeclareMathOperator{\cut}{cut}
\DeclareMathOperator{\median}{median}

\DeclareMathOperator{\NCC}{NCC}
\DeclareMathOperator{\dive}{div}
\DeclareMathOperator{\vol}{vol}

\DeclareMathOperator{\sign}{sign}

\DeclareMathOperator{\diag}{diag}
\DeclareMathOperator{\Sgn}{Sgn}

\DeclareMathOperator{\dint}{\displaystyle\int}
\newenvironment{sistema}%
{\left\{\begin{array}{@{}l@{}}}{\end{array}\right.}
\patchcmd{\abstract}{\scshape\abstractname}{\textbf{\abstractname}}{}{}
\makeatletter 
\def\@makefnmark{} 
\makeatother 

\title[Pseudo-orthogonality for $1$-Laplacian Eigenvectors and Applications ]{The Pseudo-orthogonality for Graph $1$-Laplacian Eigenvectors and Applications to Higher Cheeger Constants and Data Clustering 
}
\author[A. Corbo Esposito, G. Piscitelli]{
	Antonio Corbo Esposito$^*$, Gianpaolo Piscitelli$^*$}
\begin{document}
\maketitle
\small{
\begin{center}
{\it $^*$Dipartimento di Ingegneria Elettrica e dell'Informazione \lq\lq M. Scarano\rq\rq, \\
Universit\`a degli Studi di Cassino e del Lazio Meridionale\\ Via G. Di Biasio n. 43, 03043 Cassino (FR), Italy.}\\
\end{center}
}
\footnote{\noindent
Email: antonio.corboesposito@unicas.it, gianpaolo.piscitelli@unicas.it {\it (corresponding author)}.}

\begin{abstract} 
The data clustering problem consists in dividing a data set into prescribed groups of homogeneous data. This is a NP-hard problem that can be relaxed in the spectral graph theory, where the optimal cuts of a graph are related to the eigenvalues of graph $1$-Laplacian.
In this paper, we firstly give new notations to describe the paths, among critical eigenvectors of the graph $1$-Laplacian, realizing sets with prescribed genus.

We introduce the pseudo-orthogonality to characterize $m_3(G)$, a special eigenvalue for the graph $1$-Laplacian. Furthermore, we use it to give an upper bound for the third graph Cheeger constant $h_3(G)$, that is $h_3(G) \le m_3(G)$. This is a first step for proving that the $k$-th Cheeger constant is the minimum of the $1$-Laplacian Raylegh quotient among vectors that are pseudo-orthogonal to the vectors realizing the previous $k-1$ Cheeger constants.

Eventually, we apply these results to give a method and a numerical algorithm to compute $m_3(G)$, based on a generalized inverse power method. 


\end{abstract}

\section{Introduction}
The graph $1$-Laplacian has been deeply studied in recent years, starting from the pioneering works of Hein and B\"uhler \cite{BHa,BHb}. The study of Laplacian eigenvalues on graphs has applications in data clustering, that is the problem of dividing a data set into prescribed groups of homogeneous data. This is a NP-hard problem that can be relaxed in the spectral graph theory, where it is understood as the problem of dividing a graph into a prescribed number of groups of nodes which are densely connected inside and have little connection in between. The  clustering quality improves if we consider the eigenvalues of $p$-Laplacian as $p\to 1$, see \cite{A,BHa}. This problem has also been treated in the continuous Euclidean case \cite{BP, Ca, Che, KF, LS, Pa} and in  the anisotropic case \cite{BFK,DGPb, KN}, that is when $\R^n$ is equipped with a Finsler metric.

Furthermore, we recall that also the limiting problem of the $p$-Laplacian as $p\to\infty$ has been investigated (see \cite{JLM,EKNT} for the Euclidean case and \cite{BKJ,Pib} for the Finsler case).
Related results are obtained when other operators and boundary conditions hold (see e.g. \cite{DGPa,DP,Pia}).

Let $G=(V,E)$ an un-oriented connected planar graph, where $V$ is the vertex set, $|V|=n$, and $E\subseteq V\times V$ is the edge set. We denote by $i \sim j$ a couple of adjacent vertices $(i, j)\in E$. We study the $1$-Laplacian graph operator, that is defined as
\begin{equation}
\label{graph_1}
(\Delta_1 {\bf f})_i:=\left\{\sum_{\substack{i,j\in V\\ j\sim i}} z_{ij}({\bf f})\ \big|\  z_{i,j}({\bf f}) \in \Sgn (f_i-f_j), \ z_{ji}({\bf f })=-z_{ij}({\bf f}),\ \forall\ j\sim i\right\}, \ i=1,..,n,
\end{equation} for any ${\bf f}\in \R^n$, 
where
\begin{equation*}
\Sgn (t) =\begin{cases}
\{1\}, &t>0,\\
[-1,1],&t=0,\\
\{-1\},&t<0.\\
\end{cases}
\end{equation*}
We remark that also other definitions of graph $1$-Laplacian exist (see e.g. to \cite{A,BHa, Chu, vL} for references).

For any $i\in E$, we set
\[
d_i =\left|\{ j \in V \ | \  (i,j) \in E\} \right|, \quad i=1,...,n.
\]
Then the $1$-Laplacian eigenvalue problem is to solve a real number $\mu(G)$ and a vector ${\bf f}\in\R^n$ (respectively called eigenvalue and eigenvector associated to $\mu(G)$ of \eqref{graph_1} on $G$) satisfying 
	\begin{equation}
	\label{eigen}
	{\bf 0}\in\Delta_1 {\bf f}-\mu(G) D \Sgn ({\bf f}),
	\end{equation}
where $D:=\diag (d_1,...,d_n)$, $d:=\sum_{i\in V}d_i$ and $\Sgn({\bf f})=(\Sgn(f_1),...,\Sgn(f_n))^T$.
	
The study of eigenvectors of the $1$-Laplacian is related to the critical values of the function 
\begin{equation}
\label{I}
I({\bf f})=\sum_{\substack{i,j\in V\\ i\sim j}} |f_i-f_j|,
\end{equation}
on the symmetric piecewise linear manifold
\[
X=\left\{ {\bf f} \in \R^n \ : \ ||{\bf f}||_w=:\sum_{i=1}^n d_i |f_i|=1\right\},
\]
where $||{\bf f}||_w$ is called the $L^1-$weighted norm of ${\bf f}\in\R^n$.

In \cite{Cha, CSZZ}, the authors extended the Liusternik-Schnirelmann theory to the study of the critical points of \eqref{I}.  We us consider min-max formulas (introduced by \cite{Cu,DR}) relying on topological index theories  (see \cite{R}) involving the notion of Krasnoselskii genus. 
 
 The genus of a symmetric (i.e., $A=-A$) subset $A$ of $\R^n\setminus\{0\}$, is defined as
\begin{equation*}
 \gamma(A)=
 \begin{cases}
 0, \qquad\qquad\qquad\qquad\qquad\qquad\qquad\qquad\qquad\qquad\ \text{if}\ A=\emptyset,\\
 \min\{ k\in\N_+ \ | \ \exists\ \text{odd continuous}\  h: A \to \mathbb S^{k-1}\}, \ \text{otherwise}.
 \end{cases}
 \end{equation*}
The eigenvalues are characterized as compact paths along the spectrum $\sigma (G)$ on the symmetries of the even functional \eqref{I}. Specifically, in \cite{Cha, CSZZ}, at least $n$ critical values are obtained:
\begin{equation*}
c_k(G)=\inf_{\gamma(A)\geq k}\max_{{\bf f}\in A} \hat I({\bf f}), \quad k=1,...,n.
\end{equation*}
These eigenvalues can be ordered as
\[
c_1(G)\leq ... \leq c_n(G),
\]
but, unfortunately, it is not known if they do exhaust all the spectrum (see \cite[Sec. 6]{CSZb} for a counterexample). In this context, we denote the by $K$ the set of all critical points of $\hat I$. Then, if there exist $k,l\in \N$ such that $0\leq k\leq n-l$ and 
\[
c(G)=c_{k+1}(G)=...=c_{k+l}(G), 
\]
then $\gamma(K\cap I^{-1}(c(G))\geq l$. Furthermore, we say the eigenvalue $c(G)$ has topological multiplicity $l$ if $\gamma(K\cap \hat I^{-1}(c(G))= l$, denoted $tm(c(G))=l$.
We remark that other min-max characterization holds for eigenvalues, see e.g. \cite{DGPb} and reference therein.

Further, for any $k\ge 2$, we denote
\begin{equation}
\label{ort_value}
m_k(G)=\min_{\substack{\hat{\bf g}\in X,\ \hat{\bf g}\neq \hat{\bf g}_1,...,\hat{\bf g}_{k-1}\\ \hat{\bf g} \perp_p \hat{\bf g}_1, ... , \hat{\bf g} \perp_p\hat{\bf g}_{k-1}}}I(\hat{\bf g}),
\end{equation}
where $\hat{\bf g}_1$ is the first eigenfunction of \eqref{graph_1}, $\hat{\bf g}_j$, $j=2,...,k-1$, are inductively defined as the vectors achieving $m_j(G)$. For any $k\geq 1$, we say that $\hat{\bf g}$ and $\hat{\bf g}_1$ are pseudo-orthogonal and we denote $\hat{\bf g}\perp_p\hat {\bf g}_1$ when
\begin{equation*}
0\in \langle D \Sgn (\hat{\bf g}), \hat{\bf g}_1\rangle=\sum_{i=1}^n d_i\Sgn(\hat g_i)(\hat g_1)_i,
\end{equation*}
where $\langle \cdot,\cdot \rangle$ denotes the usual Euclidean scalar product in $\R^n$. Let us remark that the pseudo-orthogonality generalizes the condition of zero median when $k=2$. For more precise definition, see Definition \ref{ort_def} in Section \ref{aSCCC}.

The main objective of this paper is in investigating the case in which \eqref{ort_value} is equal to the $k$-Cheeger constant (particularly, the case $k=3$):
\begin{equation}
\label{cheeger_constant}
h_k(G)=\min_{S_1,S_2,...,S_k \ \text{partition of}\ V}\max_{ 1\leq i\leq k} \frac{|\partial S_i|}{\vol(S_i)}, \quad k=1,...,n.
\end{equation}
For any $A\subseteq V$, we have denoted
\[
vol(A):=\sum_{i\in A}d_i
\]
the volume of $A$ and
 \[
\partial A:=\{e=(i,j)\in E \ | \ \text{either}\  i\in A, j\notin A\ \text {or}\ j\in A, i\notin A\}
\]
the edge boundary of $A$. 

Moreover, we recall from \cite{CSZa, LGT}, the k-way Cheeger constant
\begin{equation}
\label{k-way-cheeger-constant}
\rho_k(G)=\min_{\substack{S_1,...,S_k\subset V\\S_i\cap S_j=\emptyset \ \forall i\neq j}}\max_{ 1\leq i\leq k} \frac{|\partial S_i|}{\vol(S_i)}, \quad k=1,...,n.
\end{equation}

Regarding the second Cheeger constant, it is known (see \cite{Cha,CSZa,CSZb}) that
\begin{equation}
\label{case2}
\mu_2(G)=c_2(G)=\rho_2(G)=h_2(G)=m_2(G).
\end{equation}

In this paper, we show a generalization of \eqref{case2} to case $k=3$ and give the basis for the generalization in higher cases. Particularly, we know the following inequality (refer to \cite{Cha, CSZb})
\begin{equation}
\label{muc}
\mu_k(G)\leq c_k(G)\quad\forall k\in\N.
\end{equation}
Moreover, we give a detailed proof of the following inequality (see {\bf Theorem \ref{ck_rhok}}) 
\begin{equation}
\label{cp}
c_k(G)\leq \rho_k(G)\quad\forall k\in\N.
\end{equation}
This is a known result (\cite[Th. 8]{CSZb}) but we give a proof that makes a smart use of two known results. Firstly, we generalize for the third and higher critical eigenvalues $c_k(G)$ the description in \cite{Cha} of a path joining the first and the second eigenvector to characterize sets with genus $2$. Then, we construct the paths realizing set with genus $k\ge 3$ also by using the paths joining each eigenvector with (one of) its positive part(s) as in \cite{CSZb}.

Furthermore, we recall that the reverse of \eqref{cp} holds when at least one eigenfunction associated to $c_k(G)$ has $k$ nodal domains \cite[Th. 8]{CSZb}. On the other hand, it is easily seen the following inequality
\begin{equation}
\label{rh}
\rho_k(G)\leq h_k(G)\quad\forall k\in\N.
\end{equation}
Therefore, by \eqref{muc}, \eqref{cp}, \eqref{rh}, we have
 \begin{equation*}
\mu_k(G)\leq c_k(G)\leq \rho_k(G)\le h_k(G)\quad\forall k\in\N.
\end{equation*}

Furthermore, in {\bf Theorem \ref{hNCC_thm}}, we prove 
\[
h_3(G)\leq m_3(G).
\]
Once proven this last equality, we are in position to state the main result.
 \begin{thm}
 \label{realizing} 
 Let $G=(V,E)$ be a graph, then
\begin{equation}
\label{3chain}
\mu_3(G)\leq c_3(G)\leq \rho_3(G)\leq h_3(G)\leq m_3(G).
\end{equation}
 \end{thm}

The paper is organized as follows. In the next Section, we give definitions and preliminary results on the graph $1$-Laplacian eigenvalue problem. In Section \ref{PaE}, we describe sets of prescribed genus realizing critical eigenvalues. Furthermore, in Section \ref{aSCCC} we show a suitable characterization of Cheeger constants based on pseudo-orthogonality and prove the main Theorem. Eventually, in Section \ref{AotIPMtSPD}, we show an application of these results to spectral data clustering, based on the inverse power method.

\section{The Graph 1-Laplacian Eigenvalue Problem}
Throughout this paper, for any subset $A\subseteq V$, we denote ${\bf 1}_A$ the characteristic function
\begin{equation*}
({\bf 1}_A)_i=
\begin{cases}
1, \quad i\in A,\\
0,\quad i\notin A,
\end{cases}
\end{equation*}
and $\hat {\bf 1}_A$  the normalized characteristic function
\[
\hat {\bf 1}_A= \frac{{\bf 1}_A}{\vol (A)}.\]

\begin{prop}[Cor. 2.5 \cite{Cha}]
Let $G$ be a graph and ${\bf f}$ an eigenvector associated to $\mu(G)$, then 
\[
\hat I({\bf f})=\mu(G).
\]
\end{prop}
The system \eqref{eigen} can be re-written (see \cite{Cha}) in the coordinate form as
\begin{equation}	\label{eigencoord}
\begin{cases}
\sum_{\substack{i,j\in V\\ j\sim i}} z_{ij}\in\mu(G) d_i \Sgn (f_i), \quad i=1,...,n,\\
z_{i,j}\in \Sgn (f_i-f_j),\\
z_{i,j}=-z_{i,j}.\\
\end{cases}
\end{equation}
The set of all eigenvectors, i.e., all solutions of the system \eqref{eigen} (or equivalently \eqref{eigencoord}), is denoted by $\sigma(G)$. The spectrum of a graph $G$ is finite and the elements can be ordered as
\[
\mu_1(G)\leq\mu_2(G)\leq...\leq \mu_m(G),
\]
where each eigenvalue is repeated according to its multiplicity and $m\geq n$.

\subsection{The first eigenpair}
Now we recall from \cite{Cha} some known results on eigenvalues and eigenvectors of graph $1$-Laplacian.
\begin{prop}
Let $G=(V,E)$ be a graph, then
\begin{enumerate}
\item all eigenvalues $\mu(G)$ of $\Delta_1$ satisfy $0\leq \mu(G)\leq 1$;
\item it holds $\mu(G)=1$ if and only if any nodal domain of the associated eigenvector $\bf f$ consists of a single vertex;
\item it holds $\mu(G)<1$ if and only if any nodal domain of the associated eigenvector $\bf f$ consists at least of a pair of adjacent vertices;
\item if $0<\mu(G)<1$, then $\frac 2 d\leq \mu(G)\leq \frac{n-2}{n-1}$.
\end{enumerate}
\end{prop}
We also recall that the first eigenvalue is equal to zero and the first eigenvector is constant.
\begin{prop}
Let $G=(V,E)$ be a graph, then
\begin{enumerate}
\item the first eigenvalue $\mu_1(G)=0$ is simple;
\item the eigenvector associated to $\mu_1(G)=0$ is $\hat {\bf 1}_V:=\frac 1 d {\bf 1}_V=\frac 1 d (1,...,1)$.
\end{enumerate}
\end{prop}
\begin{rem}
Let us stress that in this paper we study only connected graph. If $G$ consists of $r$ connected components, then the eigenvalue $\mu(G)=0$ has topological multiplicity $r$.
\end{rem}

\subsection{The role of the nodal domains}
To study the second and the higher eigenvalues of the graph $1$-Laplacian, it is fundamental to be able to classify the vertices of the graph in groups according the signature of any prescribed vector $\bf f$. We call nodal positive, nodal negative and null domains
\[
D_{\bf f}^+ := \{i\in V \ | \ f_i >0\}, \qquad D_{\bf f}^- := \{i\in V \ | \  f_i <0\}, \qquad D_{\bf f}^0:= \{i\in V \ | \  f_i =0\},
\]
respectively. Let $r^+({\bf f})$ and $r^-({\bf f})$ be the numbers of positive and negative nodal domains and $r({\bf f}):=r^+({\bf f})+r^-({\bf f})$. We have the following decomposition:
\[
V=\left(\bigcup_{\alpha=1}^{r^+({\bf f})} (D^+_{\bf f})_\alpha\right)\bigcup \left(\bigcup_{\beta=1}^{r^-({\bf f})}(D^-_{\bf f})_\beta\right)\bigcup D_{\bf f}^0.
\]
We denote
\[\delta_{\bf f}^\pm:=\sum_{i\in D_{\bf f}^\pm}d_i,\quad \delta_{\bf f}:=\delta_{\bf f}^++\delta_{\bf f}^-\quad\text{and}\quad\delta_{\bf f}^0:=d-\delta_{\bf f}. \] Let us observe that $d=\delta_{\bf f}^++\delta_{\bf f}^-+\delta_{\bf f}^0$. Moreover, from \cite{TH}, we recall the following nodal domain Theorem.
\begin{prop} Let $G=(V,E)$ be a graph and $\hat{\bf f}_k$ be an eigenvector associated to $c_k(G)$, $k\in\N$. If $tm(c_k)=l$, then
\[
r(\hat{\bf f}_k)\leq k+l-1.
\]
\end{prop}

\subsection{The first and the second Cheeger constants}
We recall that, in \cite[Th. 2.6]{Cha}, it has been proved that any nonconstant eigenvector of \eqref{graph_1} is pseudo-orthogonal to the first one, that means that has zero weighted median. Here we state this result for which we include the proof for the sake of completeness.
\begin{prop}
\label{second_ort} 
Let $G=(V,E)$ be a graph, $\bf f$ be an eigenvector associated to the eigenvalue $\mu (G)\neq 0$. Then ${\bf f}$ zero null weighted median: \begin{equation*}
{ 0} \in \langle D \Sgn ({\bf f}), {\bf 1}\rangle=\sum_{i=1}^nd_i\Sgn(f_i),
\end{equation*}
or, equivalently
\[
|\ \delta_{\bf f}^+-\delta_{\bf f}^-|\leq\delta_{\bf f}^0.
\]
\end{prop}
\begin{proof}
By \eqref{eigencoord}, we have
\[
\sum_{\substack{i,j\in V\\ j\sim i}} z_{ij}\in\mu(G) d_i \Sgn (f_i), \quad i=1,...,n.
\]
Therefore, since
\[
\sum_{i=1}^n\sum_{\substack{i,j\in V\\ j\sim i}} z_{ij}=0,
\]
then
\[
{ 0} \in\sum_{i=1}^nd_i\Sgn(f_i)
\]
and the conclusion follows.
\end{proof}

In this paper, we investigate the deep relationship between eigenvalues and Cheeger constants for graphs, especially related to the study of the number of nodal domains. Our aim is in generalizing to higher indices the following equality result holding for the second Cheeger constant (see \cite{Cha,CSZb}).
\begin{prop}\label{2_eq_chain}
Let $G=(V,E)$ be a connected graph, then
\[
0<c_2(G)=\mu_2(G)=\rho_2(G)=h_2(G)=m_2(G),
\]
where
\[
m_2(G)=\min_{\substack{\hat{\bf g}\in X, \hat{\bf g}\neq \hat{\bf 1}_V\\ \hat{\bf g} \perp_p \hat{\bf 1}_V}}I(\hat{\bf g}).
\]
\end{prop}
We remark that it is redundant to ask $\hat{\bf g}\neq \hat{\bf 1}_V$, because $\hat{\bf 1}_V$ is not pseudo-orthogonal to itself.

\section{Paths among Eigenvalues}\label{PaE}
In this Section, we describe how to construct paths among eigenvectors in the sublevel set of the corresponding higher eigenvalue.

\subsection{New notations to treat 1-Laplacian}
We firstly give the notations to describe the paths among eigenvalues. Since, by Proposition \ref{second_ort}, each eigenvector is equivalent to the normalized characteristic function of (one of) the positive nodal domain, we give the result for the normalized eigenvectors with only one positive nodal domain. To this aim, for each couple of vectors $\bf f$ and $\bf g$, we set
 \begin{equation}\label{degree}
 \alpha := \sum_{i\in D_{\bf f}^+\cap D_{\bf g}^+}d_i,\quad
 \beta:= \sum_{i\in D_{\bf f}^+\cap D_{\bf g}^0}d_i,\quad
 \gamma:= \sum_{i\in D_{\bf f}^0\cap D_{\bf g}^+}d_i,\quad
 \epsilon:= \sum_{i\in D_{\bf f}^0\cap D_{\bf g}^0}d_i.
 \end{equation}
They are the degrees of the intersections of the nodal domains of $\bf f$ and $\bf g$, as represented in the following table:
 \begin{center}
 \begin{tabular}{c|c|c|}
 & $D_{\bf f}^+$ & $D_{\bf f}^0$   \\ \cline {1-3}
 $D_{\bf g}^+ $ & $\alpha $& $\gamma$ \\ \cline {1-3}
 $D_{\bf g}^0 $ & $\beta$ & $\epsilon $ \\ \cline {1-3}
 \end{tabular}
 \end{center}
Moreover, we denote by ${\bf E^+_{\bf f}}=\sum_{i\in D_{\bf f}^+}{\bf e}_i$, ${\bf E}^0_{\bf f}=\sum_{i\in D_{\bf f}^0}{\bf e}_i$, ${\bf E}_\alpha=\sum_{i\in D_{\bf f}^+\cap D_{\bf g}^+}{\bf e}_i$, ${\bf E}_\beta=\sum_{i\in D_{\bf f}^+\cap D_{\bf g}^0}{\bf e}_i$, ${\bf E}_\gamma=\sum_{i\in D_{\bf f}^0\cap D_{\bf g}^+}{\bf e}_i$, ${\bf E}_\epsilon=\sum_{i\in D_{\bf f}^0\cap D_{\bf g}^0}{\bf e}_i$.

Furthermore, we partitionate the set of edges $E$ in ten subsets. We denote the subset of couples of $E$ for which the indices are in the same intersection of nodal domains as
 \begin{equation}\label{edge_same}
 \begin{split}
A:=\{e=(i,j)\in E\ | \ i,\ j \in D_{{\bf f}}^+\cap D_{{\bf g}}^+ \},\\
B:=\{e=(i,j)\in  E \ | \  i,\ j\in D_{{\bf f}}^+\ \cap D_{{\bf g}}^0\},\\
C:=\{e=(i,j)\in E \ | \   i,\ j\in D_{{\bf f}}^0\ \cap  D_{{\bf g}}^+\},\\
D:=\{e=(i,j)\in E \ | \   i,\ j\in D_{{\bf f}}^0\ \cap  D_{{\bf g}}^0\}.\\
 \end{split}
 \end{equation}
The subset of couples of $E$ for which the indices are in one nodal domain  but in different intersections of nodal domains are denoted as
  \begin{equation}\label{edge_diff}
 \begin{split}
\tilde E:=\{e=(i,j)\in D_{{\bf f}}^+\times D_{{\bf f}}^+\subset E \ | \ \text{either}\  i\in D_{{\bf g}}^+,\ j\in D_{{\bf g}}^0\ \text {or}\ j\in D_{{\bf g}}^0,\ i\in D_{{\bf g}}^+\},\\
F:=\{e=(i,j)\in D_{{\bf f}}^0\times D_{{\bf f}}^0\subset E \ | \ \text{either}\  i\in D_{{\bf g}}^+,\ j\in D_{{\bf g}}^0\ \text {or}\ j\in D_{{\bf g}}^0,\ i\in D_{{\bf g}}^+\},\\
G:=\{e=(i,j)\in D_{{\bf g}}^+\times D_{{\bf g}}^+\subset E \ | \ \text{either}\  i\in D_{{\bf f}}^+,\ j\in D_{{\bf f}}^0\ \text {or}\ j\in D_{{\bf f}}^0,\ i\in D_{{\bf f}}^+\},\\
H:=\{e=(i,j)\in D_{{\bf g}}^0\times D_{{\bf g}}^0\subset E \ | \ \text{either}\  i\in D_{{\bf f}}^+,\ j\in D_{{\bf f}}^0\ \text {or}\ j\in D_{{\bf f}}^0,\ i\in D_{{\bf f}}^+\},\\
 \end{split}
 \end{equation}
where we use $\tilde E$ instead of $E$, since this yet denote the set of the edges. Furthermore, the subset of couples of $E$ for which both the indices are in different nodal domains and in a different intersections of nodal domains are denoted as
  \begin{equation}
 \begin{split}\label{edge_diff_diff}
L:=\{ &e=(i,j)\in E \ | \\ & \text{either}\  i\in D_{{\bf f}}^+\cap D_{{\bf g}}^+,\ j\in D_{{\bf f}}^0\cap D_{{\bf g}}^0\ \text {or}\ j\in D_{{\bf f}}^0\cap D_{{\bf g}}^0,\ i\in D_{{\bf f}}^+\cap D_{{\bf g}}^+\},\\
M:=\{ & e=(i,j)\in E \ | \\ &\text{either}\  i\in D_{{\bf f}}^0\cap D_{{\bf g}}^+,\ j\in D_{{\bf f}}^+\cap D_{{\bf g}}^0\ \text {or}\ j\in D_{{\bf f}}^+\cap D_{{\bf g}}^0,\ i\in D_{{\bf f}}^0\cap D_{{\bf g}}^+\}.\\
 \end{split}
 \end{equation}
These ten subsets of $E$ are represented in the following table
 
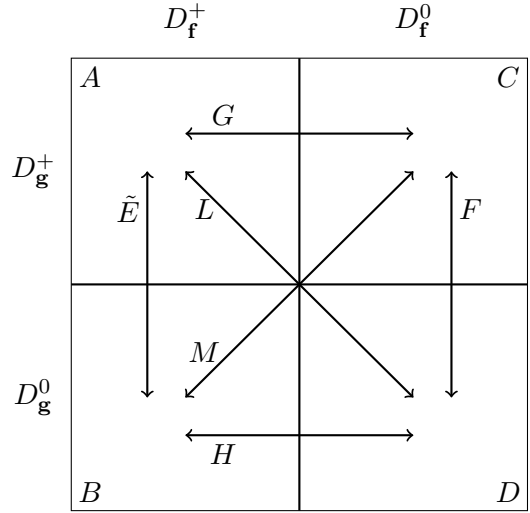
\begin{center}
 \begin{tikzpicture}
\draw (-3,-3) rectangle (3,3);
\draw[thick,-] (-3,0) -- (3,0);
\draw[thick,-] (0,-3) -- (0,3);
\draw[thick,<->] (-2,-1.5) -- (-2,1.5);
\draw[thick,<->] (2,-1.5) -- (2,1.5);
\draw[thick,<->] (-1.5,2) -- (1.5,2);
\draw[thick,<->] (-1.5,-2) -- (1.5,-2);
\draw[thick,<->] (-1.5,-1.5) -- (1.5,1.5);
\draw[thick,<->] (-1.5,1.5) -- (1.5,-1.5);
\node[align=left] (punto) at (-2.75,2.75) {$A$};
\node[align=left] (punto) at (-2.75,-2.75) {$B$};
\node[align=left] (punto) at (2.75,2.75) {$C$};
\node[align=left] (punto) at (2.75,-2.75) {$D$};
\node[align=left] (punto) at (-2.25,1) {$\tilde E$};
\node[align=left] (punto) at (2.25,1) {$F$};
\node[align=left] (punto) at (-1,2.25) {$G$};
\node[align=left] (punto) at (-1,-2.25) {$H$};
\node[align=left] (punto) at (-1.25,1) {$L$};
\node[align=left] (punto) at (-1.25,-0.9) {$M$};
\node[align=left] (punto) at (-3.5,1.5) {$D_{\bf g}^+$};
\node[align=left] (punto) at (-3.5,-1.5) {$D_{\bf g}^0$};
\node[align=left] (punto) at (-1.5,3.5) {$D_{\bf f}^+ $};
\node[align=left] (punto) at (1.5,3.5) {$D_{\bf f}^0$};
\end{tikzpicture}
 \end{center}
Now, we denote $a=|A|$, $b=|B|$, $c=|C|$, $d=|D|$, $\tilde e=|\tilde E|$, $f=|F|$, $g=|G|$, $h=|H|$, $l=|L|$, $m=|M|$. By the use of this notation, the following equalities hold:
\begin{equation*}
\begin{split}
\delta_{\bf f}^+=\alpha+\beta=2a+2b+2\tilde e+g+h+l+m,\\
\delta_{\bf f}^0=\gamma+\epsilon=2c+2d+2f +g+h+l+m,\\
\delta_{\bf g}^+=\alpha+\gamma=2a+2c+2g+ \tilde e+f+l+m,\\
\delta_{\bf g}^0=\beta+\epsilon=2b+2d+2h +\tilde e+f+l+m.
\end{split}
\end{equation*}
Furthermore, if ${\bf f}$ and ${\bf g}$ are also eigenvectors, then, by Proposition \ref{second_ort}, they have zero weighted median. This property leads to the following inequalities:
\begin{equation*}
\begin{split}
& a+b+\tilde e\leq c+d+f,\\
& a+c+g\leq b+d+h.
\end{split}
\end{equation*}
By summing these inequalities we also have
\begin{equation*}
\begin{split}
& \alpha\leq\epsilon,\\
& 2a+\tilde e+g\leq 2d +f+h.
\end{split}
\end{equation*}

Moreover, if the eigenvalues $\mu_{\bf f}(G)$ and $\mu_{\bf g}(G)$ are associated to ${\bf f}$ and ${\bf g}$, respectively, we have
\begin{equation*}
\begin{split}
   \mu_{\bf f}(G)=\frac{2g+2h+2l+2m}{2a+2b+2\tilde e+g+h+l+m},\\
 \mu_{\bf g}(G)=\frac{2\tilde e+2f+2l+2m}{2a+2c+2g+\tilde e+f+l+m}.\\
\end{split}
\end{equation*}
 Using the notation in \eqref{edge_same}-\eqref{edge_diff}-\eqref{edge_diff_diff}, we have that
 \begin{equation*}
 \begin{split}
 \delta^+_{\bf f}= 2 a+2b+2\tilde e+g+h+l+m,\\
  \delta^+_{\bf g}= 2 a+2c+2g+\tilde e+f+l+m.
 \end{split}
 \end{equation*}
 
\subsection{The behaviour of the eigenvectors}
\begin{prop} 
Let $G$ be a graph and ${\bf f}=\hat{\bf 1}_{D_{{\bf f}}^+}$ and ${\bf g}=\hat{\bf 1}_{D_{{\bf g}}^+}$ two eigenvectors, respectively associated to $\mu_{\bf f}(G)$ and $\mu_{\bf g}(G)$. Then
   \begin{align}
   \label{destra_f} &{ 0}\in  \langle D \Sgn ({\bf f}), {\bf g}  \rangle - \langle D \Sgn( {\bf g}), {\bf f}  \rangle;\\ 
   \label{sinistra_f} &0\in  \langle \Delta_1 {\bf f},  {\bf g} \rangle - \langle \Delta_1  {\bf g}, {\bf f}  \rangle;\\ 
\label{diff_f_dim1}& 0\in\langle\Delta_1 {\bf f},{\bf g}\rangle-\mu_{\bf f}(G) \langle D \Sgn ({\bf f}),{\bf g}\rangle;\\  
\label{diff_f_dim2}&  0\in\langle\Delta_1 {\bf g},{\bf f}\rangle-\mu_{\bf g}(G) \langle D \Sgn ({\bf g}),{\bf f}\rangle.
 \end{align} 
 \end{prop}
 \begin{proof}  Using the notation in \eqref{degree}, we have that
\begin{equation*}
   \langle D \Sgn ({\bf f}), {\bf g}  \rangle-   \langle D \Sgn( {\bf g}), {\bf f}  \rangle=\frac {\alpha(\beta-\gamma)+\gamma(\alpha+\beta)\Sgn (0)+\beta(\alpha+\gamma) \Sgn (0)} {(\alpha+\gamma)(\alpha+\beta)}
\end{equation*}
contains $0$ and \eqref{destra_f} is proved.

 Using the notation in \eqref{edge_same}-\eqref{edge_diff}-\eqref{edge_diff_diff}, we have that
 \begin{equation*}
 \begin{split}
 \delta^+_{\bf f}= 2 a+2b+2\tilde e+g+h+l+m,\\
  \delta^+_{\bf g}= 2 a+2c+2g+\tilde e+f+l+m.
 \end{split}
 \end{equation*}
We have
 \begin{equation*}
 \begin{split}
  \langle \Delta_1 {\bf f},  {\bf g} \rangle =\sum_{i\in V} \sum_{\substack{j\in V\\ i\sim j}}  \Sgn (f_i-f_j)\ g_i= \frac{(l-m)+(2a+2c+\tilde e+f)\Sgn (0)}{\delta^+_{\bf g}},\\
   \langle \Delta_1  {\bf g}, {\bf f}  \rangle= \sum_{i\in V}\sum_{\substack{j\in V\\ i\sim j}} \Sgn (g_i-g_j)\ f_i= \frac{(l-m)+(2a+2b+g+h)\Sgn (0)}{\delta^+_{\bf f}}.
 \end{split}
 \end{equation*}
 Hence the difference  
 \begin{equation*}
 \begin{split}
  \langle \Delta_1 {\bf f},  {\bf g} \rangle -&   \langle \Delta_1  {\bf g}, {\bf f}  \rangle=  \frac{(l-m)(2b+\tilde e+h-2c-g-f)}{\delta^+_{\bf f}\delta^+_{\bf g}}\\
  &+\frac{[(2a+2c+e+f)(2 a+2b+2\tilde e+g+h+l+m)\Sgn (0)}{\delta^+_{\bf f}\delta^+_{\bf g}}\\
  &+\frac{(2a+2b+g+h)(2 a+2c+2g+\tilde e+f+l+m)\Sgn (0)}{\delta^+_{\bf f}\delta^+_{\bf g}}
 \end{split}
 \end{equation*}
 contains $0$ and \eqref{sinistra_f} follows.

The eigenvectors  ${\bf f}$ and ${\bf g}$ satisfies
\begin{equation*}
	{\bf 0}\in\Delta_1 {\bf f}-\mu_{\bf f}(G) D \Sgn ({\bf f}) \quad\text{and}\quad
	{\bf 0}\in\Delta_1 {\bf g}-\mu_{\bf g}(G) D \Sgn ({\bf g}). 
\end{equation*}
 By multiplying the first relation by ${\bf g}$ and the second one by ${\bf f}$, we have \eqref{diff_f_dim1}-\eqref{diff_f_dim2}.
   \end{proof}

 \subsection{Paths among eigenvalues}
To characterize the sets realizing the third (and higher) critical eigenvalues, we need to construct paths between eigenvectors. 
In this context, for any $c\in\R $, we denote the level and the sublevel set of $I$ on $X$ as
\[
\hat I_c=\{ {\bf f}\in X \ |\ I( {\bf f} )=c \}\quad\text{and}\quad\hat I^-_c=\{ {\bf f}\in X \ |\ I( {\bf f} )\le c \},
\]
respectively.

\begin{defn} 
Let $G=(V,E)$ be a graph. For any $A\subseteq V$, we say that ${\bf f}$ and ${\bf h}\in\R^n$ are equivalent in A and we denote ${\bf f}\simeq {\bf h}$ in $A$, if there exists a path $\gamma(t)$ in $X$ such that $\gamma(0)={\bf f}$, $\gamma(1)={\bf h}$ and $\gamma(t)\in A$ for any $t\in[0,1]$, . 
\end{defn}
\begin{prop}[Th.1 \cite{CSZb}]
\label{percorsi_pos} 
Let $G=(V,E)$ be a graph, $\bf f$ be an eigenvector of \eqref{graph_1} associated to the eigenvalue $\mu_{\bf f} (G)\neq 0$. Then the positive (or negative) part of an eigenvector realizes the same eigenvalue, that is $\forall\ \alpha \in (\{1,...,r^+({\bf f})\})$ (or $\forall\ \beta\in \{1,...,r^-({\bf f})\})$, we have
 \[
 {\bf f} \simeq {\bf \hat 1}_{D_\alpha^+} \quad \text{(or} \ {\bf f} \simeq {\bf \hat 1}_{D_\beta^-})\quad \text{in} \ \sigma(G)\cap\hat I^-_{\mu_{\bf f} (G)}.
 \]
\end{prop} 
  \begin{prop}
  \label{percorsi_3} Let $G=(V,E)$ be a graph and $\hat {\bf 1}_V$, $\hat{\bf f}_2$ and $\hat{\bf f}_3$ eigenvectors of \eqref{graph_1} associated to the eigenvalues $\mu_1(G)$, $\mu_2(G)$ and $\mu_3(G)$, respectively. Then
 \begin{enumerate}
 \item $\hat{\bf 1}_V\simeq\hat {\bf f}_2$ in $\hat I^-_{\mu_2(G)}$; 
 \item $\hat{\bf 1}_V\simeq\hat {\bf f}_3$ in $\hat I^-_{\mu_3(G)}$; 
 \item $\hat{\bf f}_2\simeq\hat {\bf f}_3$ in $\hat I^-_{\mu_3(G)}$. 
 \end{enumerate}
 \end{prop}
 \begin{proof} 
To prove item {\it (1)} we remark that, by Proposition \ref{percorsi_pos}, $\hat{\bf f}_2\simeq \hat{\bf 1}_{D_{{\bf f}_2}^+}$ in $\hat I_{\mu_2}^-$, then it is sufficient to prove that $\hat{\bf 1}_V\simeq \hat{\bf 1}_{D_{\hat{\bf f}_2}^+}$ in $\hat I_{\mu_2}^-$. Let us consider 
 \[
 \varphi(t)=
 t \hat{\bf 1}_{D_{{\bf f}_2}^+} +(1-t)\hat{\bf 1}_V=\left(\frac t{\delta^+_{{\bf f}_2}} +\frac{1-t}d\right){\bf E}^+_{{\bf f}_2}+\frac{1-t}d{\bf E}^0_{{\bf f}_2}
 , \quad t\in[0,1].
 \] 
 We set 
 \begin{equation*}
 \begin{split}
E_1=\{e=(i,j)\in E \ | \ \text{either}\  i\in D_{{\bf f}_2}^+,\ j\in D_{{\bf f}_2}^0\ \text {or}\ i\in D_{{\bf f}_2}^0,\ j\in D_{{\bf f}_2}^+\}.\\
 \end{split}
 \end{equation*}
 Then, we have $I (\varphi (t))=\frac{t |E_1|}{\delta^+_{{\bf f}_2}
 }$ and $||t \hat{\bf 1}_{D_{{\bf f}_2}^+} +(1-t)\hat{\bf 1}_V||_w=1$ and hence
 \[
\hat I (\varphi (t))=\frac{t |E_1|}{\delta^+_{{\bf f}_2}}.
 \]
Therefore the conclusion follows by noting that \[
I (\varphi (t))'=\frac{|E_1|
}{\delta^+_{{\bf f}_2}
}>0,
\]
and that $$I(\varphi(0))=I(\hat{\bf 1}_V)=\mu_1(G)=0\leq \mu_2(G)=I( \hat{\bf f}_2)=I(\varphi (1)).$$ The proof of item {\it (2)} follows analogously. To prove item {\it (3)}, again by Proposition \ref{percorsi_pos}, it is sufficient to prove that $\hat{\bf 1}_{D_{{\bf f}_2}^+}\simeq \hat{\bf 1}_{D_{{\bf f}_3}^+}$ in $\hat I_{\mu_3}^-$.
Let us consider the following path
 \[
\psi(t)={t \hat{\bf 1}_{D_{{\bf f}_3}^+} +(1-t)\hat{\bf 1}_{D_{{\bf f}_2}^+}}=\left(\frac{1-t}{\delta^+_{{\bf f}_2}}+\frac t {\delta^+_{{\bf f}_3}}\right){\bf E}_\alpha+\frac t {\delta^+_{{\bf f}_3}}{\bf E}_\beta+\frac{1-t}{\delta^+_{{\bf f}_2}}{\bf E}_\gamma, \quad t\in[0,1].\] 
Using the notation in \eqref{degree} and \eqref{edge_same}-\eqref{edge_diff}-\eqref{edge_diff_diff}, with ${\bf g}=\hat{\bf f}_2$ and ${\bf f}=\hat{\bf f}_3$, we have
 \begin{equation}\label{mu2mu3}
 \begin{split}
 \frac{2\tilde e+2f+2l+2m}{\delta^+_{\hat{\bf f}_2}}&=I(\psi (0))=I(\hat{\bf f}_2) =\mu_2(G)\le\\
& \leq  \mu_3(G)=I(\hat{\bf f}_3)=I(\psi(1))=\frac{2g+2h+2l+2m}{\delta^+_{{\bf f}_3}}.
 \end{split}
 \end{equation}
Then, we have 
\begin{equation*}
I (\psi(t))=
\begin{sistema}
\frac{2\tilde e+2f+2l+2m}{\delta^+_{{\bf f}_2}}+t\frac{\delta^+_{{\bf f}_2}(2g+2h+2l-2m)-\delta^+_{{\bf f}_3}(2 \tilde e+2f+2l+2m)}{\delta^+_{{\bf f}_2}\delta^+_{{\bf f}_3}}\quad t<\frac{  \delta^+_{{\bf f}_3}  }{    \delta^+_{{\bf f}_3}    + \delta^+_{{\bf f}_3}  },\\
\frac{2\tilde e+2f+2l-2m}{\delta^+_{{\bf f}_2}}+t\frac{\delta^+_{{\bf f}_2}(2g+2h+2l+2m)-\delta^+_{{\bf f}_3}(2 \tilde e+2f+2l-2m)}{\delta^+_{{\bf f}_2}\delta^+_{{\bf f}_3}} \quad t\ge\frac{  \delta^+_{{\bf f}_3}  }{    \delta^+_{{\bf f}_3}    + \delta^+_{{\bf f}_3}  },\\\end{sistema}
\end{equation*}
and $||t \hat{\bf 1}_{D_{{\bf f}_2}^+} +(1-t)\hat{\bf 1}_{D_{{\bf f}_3}^+}||_w=1$. Therefore $I(\psi(t))= \hat I(\psi(t))$ and 
\begin{equation*}
I' (\psi(t))=
\begin{sistema}
\frac{\delta^+_{{\bf f}_2}(2g+2h+2l-2m)-\delta^+_{{\bf f}_3}(2 \tilde e+2f+2l+2m)}{\delta^+_{{\bf f}_2}\delta^+_{{\bf f}_3}}\quad t<\frac{  \delta^+_{{\bf f}_3}  }{    \delta^+_{{\bf f}_3}    + \delta^+_{{\bf f}_3}  },\\
\frac{\delta^+_{{\bf f}_2}(2g+2h+2l+2m)-\delta^+_{{\bf f}_3}(2 \tilde e+2f+2l-2m)}{\delta^+_{{\bf f}_2}\delta^+_{{\bf f}_3}} \quad t\ge\frac{  \delta^+_{{\bf f}_3}  }{    \delta^+_{{\bf f}_3}    + \delta^+_{{\bf f}_3}  }.\\
\end{sistema}
\end{equation*}
Finally, by \eqref{mu2mu3}, the term in the second line is nonnegative and hence the conclusion follows.
 \end{proof}
We remark that the first item of this Lemma has been proven in \cite[Lem. 5.1]{Cha} by using different notations. 

Now, we can generalize the result of the previous Proposition. The proof can be easily given by following line by line the proof of Proposition \ref{percorsi_3}.
\begin{prop}\label{percorsi_k}
Let $G$ be a graph and $\hat{\bf f}_h$ and $\hat{\bf f}_k$ associated, respectively, to the eigenvalue $\mu_h(G)\leq\mu_k(G)$. Then
\[
\hat{\bf f}_h\simeq\hat {\bf f}_k \ \text{in}\ \hat I^-_{\mu_k}.
\]
\end{prop}

Furthermore, we recall the following inequality results from \cite{CSZb}. 
\begin{prop}
Let $G=(V,E)$ be a graph. 
Then
\begin{enumerate}
\item $c_k(G)\leq \rho_k(G)$, for all $k\in\{1,...,n\}$;
\item if $\hat {\bf f}_k$ is an eigenvector associated to $c_k(G)$ such that $r(\hat{\bf f}_k)\geq s$, then $\rho_s(G)\leq c_k(G)$.
\end{enumerate}
\end{prop} 

Since at this point, we are able to construct (and explicitly describe) sets with genus $k$ realizing the $k$-th critical value $c_k(G)$, we can give a detailed proof of the inequality $c_k(G)\leq \rho_k(G)$.
\begin{thm}\label{ck_rhok}
For any $k\in\N$, we have
\begin{equation*}
\mu_k(G)\leq c_k(G)\leq \rho_k(G)\le h_k(G).
\end{equation*}
\end{thm}
\begin{proof}
The first and the last inequalities are easily seen since the critical values do not exahust all the spectrum and since the class of the $k$-partition of $V$ is contained in the class of the all $k$-tuple of disjoint sets, respectively. Hence we only need to prove $c_k(G)\leq \rho_k(G)$. 
Firstly we analyze the inequality for $k=3$: $c_3(G)\leq \rho_3(G)$.

We observe that $\rho_3(G)> c_2 (G)$ and that, since $I$ is continuous and we are considering a compact set, the minimum is achieved. Hence it exists $\Psi\in\R^n$ such that $\rho_3(G)=\hat I (\Psi)$.

Therefore, even if we are not able to say that $\rho_3(G)=\mu_j(G)$ for some $j\geq 3$, we can suppose that $\Psi$ is a characteristic function of a certain domain, because
the vectors realizing the Cheeger constants are normalized characteristic functions of a certain domain $A$, indeed we have $ \frac{|\partial A|}{\vol(A)}=\hat I({\bf 1}_A)=I(\hat{\bf 1}_A)$.

Let us denote by $\hat{\bf f }_2$ and $\hat{\bf f}_1$ the (normalized positive) eigenvectors associated to $\mu_2(G)$, $\mu_1(G)=0$, respectively. By Lemma \ref{percorsi_k}, it is possible to construct a path $\gamma_1$ connecting $\Psi$ and $\hat{\bf f}_1$, a path $\gamma_2$ connecting $\hat{\bf f}_2$ and $\hat{\bf 1}_V$ and a path $\gamma_3$ connecting $\Psi$ and  $\hat{\bf f}_2$. Since it is seen that $\hat{\bf f}_1$, $\hat{\bf f}_2$ and $\Psi$ are equivalent in $\hat I^-_{\rho_3(G)}$, we consider the linear convex combination $T_1$ of vertices $\Psi$, $\hat{\bf f}_2$ and $\hat{\bf f}_1$. We have 
 \begin{equation*}
 T_1(t_1,t_2)=t_1\hat{\bf f}_1+t_2\hat {\bf f}_2+(1-t_1-t_2)\Psi,
 \end{equation*}
 with $t_1, t_2\ge 0$ and ${t_1}+t_2\leq1$ (see the graph below). Then, we have 
\begin{equation*}
I (T_1(t_1,t_2))=
\begin{sistema}
\frac{2g+2h+2l+2m}{\delta^+_{{\bf f}_3}}-t_1\frac{2g+2h+2l+2m}{\delta^+_{{\bf f}_3}}+t_2\frac{\delta^+_{{\bf f}_2}(2g+2h+2l-2m)-\delta^+_{{\bf f}_3}(2 \tilde e+2f+2l+2m)]}{\delta^+_{{\bf f}_2}\delta^+_{{\bf f}_3}}\\
\qquad\qquad\qquad\qquad\qquad\qquad\qquad\qquad\qquad\qquad\qquad t_2<\frac{  \delta^+_{{\bf f}_2} (1-t_1) }{    \delta^+_{{\bf f}_3}    + \delta^+_{{\bf f}_3}  },\\
\frac{2g+2h+2l-2m}{\delta^+_{{\bf f}_3}}-t_1\frac{2g+2h+2l-2m}{\delta^+_{{\bf f}_3}}+t_2\frac{\delta^+_{{\bf f}_2}(2g+2h+2l+2m)-\delta^+_{{\bf f}_3}(2 \tilde e+2f+2l-2m)]}{\delta^+_{{\bf f}_2}\delta^+_{{\bf f}_3}}\\
\qquad\qquad\qquad\qquad\qquad\qquad\qquad\qquad\qquad\qquad\qquad t_2\geq\frac{  \delta^+_{{\bf f}_2}(1-t_1)  }{    \delta^+_{{\bf f}_3}    + \delta^+_{{\bf f}_3}  },
\end{sistema}
\end{equation*}
and
\begin{equation}\label{one}
||T_1(t_1,t_2)||_w=1.
\end{equation}
Therefore $I(T_1(t_1,t_2))= \hat I(T_1(t_1,t_2))$ and, when  $t_2\geq\frac{  \delta^+_{{\bf f}_2}(1-t_1)  }{    \delta^+_{{\bf f}_3}    + \delta^+_{{\bf f}_3}  } $, we have
\begin{equation*}
\begin{split}
\frac{d }{dt_1}I (T_1(t_1,t_2)) & =-\frac{2g+2h+2l+2m}{\delta^+_{{\bf f}_3}}<0,\\
\frac{d }{dt_2}I (T_1(t_1,t_2)) & =\frac{\delta^+_{{\bf f}_2}(2g+2h+2l+2m)-\delta^+_{{\bf f}_3}(2 \tilde e+2f+2l-2m)}{\delta^+_{{\bf f}_2}\delta^+_{{\bf f}_3}}<0.
\end{split}
\end{equation*}
 Hence, by noting that $T_1(t_1,t_2)=t_1 \hat{\bf f}_1+t_2\hat{\bf f}_2\in  \hat I^-_{\mu_2(G)}$ when $t_1+t_2=1$, that $T_1(1,0)=\hat{\bf f}_1$, $T_1(0,1)=\hat{\bf f}_2$ and $T_1(0,0)=\Psi$, we have $I(T_1(t_1,t_2))\leq{\rho_3(G)}$, for $t_1, t_2\ge 0$ and ${t_1}+t_2\leq1$ (see the graph below). 
  
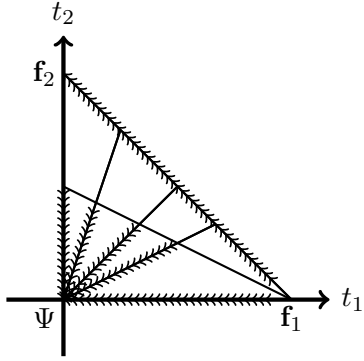
\begin{center}
 \begin{tikzpicture}
\draw[->,ultra thick] (-0.75,0)--(3.5,0) node[right]{$t_1$};
\draw[->,ultra thick] (0,-0.75)--(0,3.5) node[above]{$t_2$};
\draw[thick,-] (3,0) -- (0,3);
\draw[thick,-] (3,0) -- (0,1.5);
\draw[thick,-] (0,0) -- (2,1);
\draw[thick,-] (0,0) -- (1.5,1.5);
\draw[thick,-] (0,0) -- (0.75,2.25);
\node[align=left] (punto) at (3,-0.25) {${\bf f}_1$};
\node[align=left] (punto) at (-0.25,3) {${\bf f}_2$};
\node[align=left] (punto) at (-0.25,-0.25) {$\Psi$};
\draw[thick,->>>>>>>>>>>>>>>>>>>>>>>>>>>>>>>>>>>>>>>>>>>] (3,0) -- (0,3);
\draw[thick,->>>>>>>>>>>>>>>>>>>>>>>>>>>>>>] (3,0) -- (0,0);
\draw[thick,->>>>>>>>>>>>>>>>] (0,3) -- (0,0);
\draw[thick,<<<<<<<<<<<<<<<<<<-] (0,0) -- (2,1);
\draw[thick,<<<<<<<<<<<<<<<-] (0,0) -- (1.5,1.5);
\draw[thick,<<<<<<<<<<<<<<-] (0,0) -- (0.75,2.25);
\end{tikzpicture}
 \end{center}
 Therefore $T_1(t_1,t_2))\in \hat I^-_{\rho_3(G)}$, for $t_1, t_2\ge 0$ and ${t_1}+t_2\leq1$.

Similarly we can construct the other seven linear convex combinations $T_i$, $i=2,...,8$ with the first vertex between $\pm\Psi$, the second one between $\pm{\bf f}_2$ and the third one between $\pm{\bf f}_1$. Since the norm of $T_i$, $i=1,...,8$ is unitary as in \eqref{one}, then this result can be also showed by using the convexity of $I$. For the convenience of the reader, we preferred to give the proof without using this property.

We remark that, when considering $\pm{\bf f}_1$, $\pm{\bf f}_2$ and $\pm\Psi$, the previous construction could not give a set with genus 3. Hence, alternatively, we could also consider the other two normalized characteristic functions $\hat{1}_R$ and $\hat{1}_S$ of the triple realizing $\rho_3(G)$. Indeed, $\pm\Psi$, $\pm\hat{1}_R$ and $\pm\hat{1}_S$ are in $\hat I^-_{\rho_3(G)}$ and the span $<\Psi, \hat{1}_R, \hat{1}_S>$ has genus 3.

By gluing the eight convex combinations $T_i$, $i=1,...,8$, we obtain a set $A$ isomorphic to $\mathbb S^2$. Hence $\gamma(A)\geq 3$ and then 
 \[
 c_3(G)=\inf_{\gamma(A)\geq 3}\max_{{\bf f}\in A} I({\bf f})\leq \rho_3(G).
 \]
 Therefore the conclusion follows.
 
Finally, this result also lead to the description of sets with genus $k$ for any $k\in\N$. We consider a linear convex combination $\sum_{i=1}^Nt_i{\bf f}_i$, for $t_i\geq 0$, $\sum_{i=1}^Nt_i$, $N\in\N$. It has unitary norm, indeed 
\begin{equation*}
\left|\left|\sum_{i=1}^Nt_i{\bf f}_i\right|\right|_w=\sum_{i=1}^Nt_i\left[\sum_{
k_j\neq k_i}
\frac{|D^{k_1,...,k_{j-1},+,k_{j+1},...,k_N}|}{\delta^+_{{\bf f}_i}}\right]=\sum_{i=1}^Nt_i=1,
\end{equation*}
where $D^{k_1,...,k_{j-1},+,k_{j+1},...,k_N}=D^{k_1}_{{\bf f}_1}\cap ...\cap D^{k_{j-1}}_{{\bf f}_{{j-1}}}\cap D^{+}_{{\bf f}_{j}}\cap D^{k_{j+1}}_{{\bf f}_{{j+1}}}\cap ...\cap D_{{\bf f}_N}^{k_N}$, for $k_j\in \{+,0\}$.
 By using the convexity of $I$ and constructing analogously the iper-surfaces in any $2^n$-ant, the conclusion follows.
  \end{proof}
 We conclude this Section by analyzing the pseudo-orthogonality of eigenvectors of two special graphs: the path graph $P_{10}$ and the cycle graph $C_{10}$ in dimension $n=10$.
 \begin{exa}
 The degree of each vertices of the graph $P_{10}$ is equal to two, except for the first and the last one, for which the degree is equal to one. The spectrum is $\sigma(P_{10})=\left\{0,\frac 19,\frac 17,\frac 15,\frac 14,\frac 13, \frac 12,1\right\}$. We draw the eigenvectors in a form representing the maximum number of nodal domains. 
  
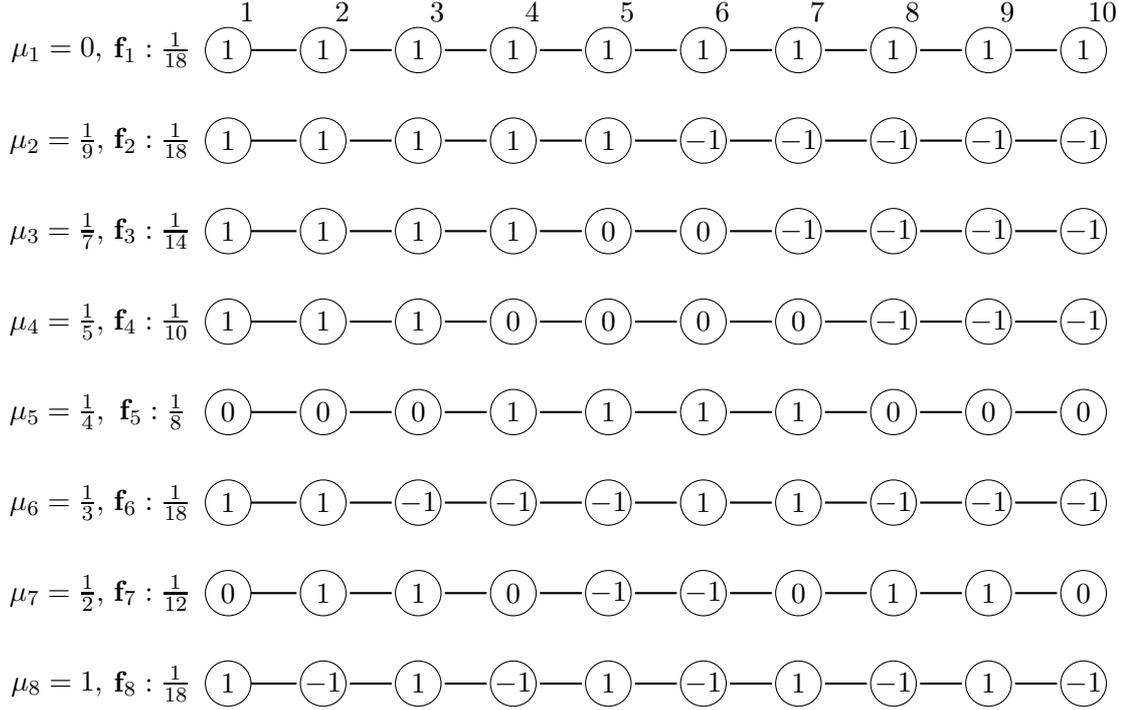
\begin{center}
\begin{tikzpicture}
\draw (-6.25,0) circle (0.3cm);
\draw (-5,0) circle (0.3cm);
\draw (-3.75,0) circle (0.3cm);
\draw (-2.5,0) circle (0.3cm);
\draw (-1.25,0) circle (0.3cm);
\draw (0,0) circle (0.3cm);
\draw (1.25,0) circle (0.3cm);
\draw (2.5,0) circle (0.3cm);
\draw (3.75,0) circle (0.3cm);
\draw (5,0) circle (0.3cm);
\draw[thick,-] (-5.95,0) -- (-5.35,0);
\draw[thick,-] (-4.65,0) -- (-4.1,0);
\draw[thick,-] (-3.4,0) -- (-2.85,0);
\draw[thick,-] (-2.15,0) -- (-1.6,0);
\draw[thick,-] (-0.9,0) -- (-0.35,0);
\draw[thick,-] (0.35,0) -- (0.9,0);
\draw[thick,-] (1.6,0) -- (2.15,0);
\draw[thick,-] (2.85,0) -- (3.4,0);
\draw[thick,-] (4.1,0) -- (4.65,0);
\node[align=left] (punto) at (-6,0.5) {$1$};
\node[align=left] (punto) at (-4.75,0.5) {$2$};
\node[align=left] (punto) at (-3.5,0.5) {$3$};
\node[align=left] (punto) at (-2.25,0.5) {$4$};
\node[align=left] (punto) at (-1,0.5) {$5$};
\node[align=left] (punto) at (0.25,0.5) {$6$};
\node[align=left] (punto) at (1.5,0.5) {$7$};
\node[align=left] (punto) at (2.75,0.5) {$8$};
\node[align=left] (punto) at (4,0.5) {$9$};
\node[align=left] (punto) at (5.25,0.5) {$10$};
\node[align=left] (punto) at (-8.5,-0) {${\mu}_1=0,$};
\node[align=left] (punto) at (-7.25,-0) {${\bf f}_1:\frac1{18}$};
\node[align=left] (punto) at (-6.25,-0) {$1$};
\node[align=left] (punto) at (-5,-0) {$1$};
\node[align=left] (punto) at (-3.75,-0) {$1$};
\node[align=left] (punto) at (-2.5,0) {$1$};
\node[align=left] (punto) at (-1.25,0) {$1$};
\node[align=left] (punto) at (0,0) {$1$};
\node[align=left] (punto) at (1.25,0) {$1$};
\node[align=left] (punto) at (2.5,-0) {$1$};
\node[align=left] (punto) at (3.75,-0) {$1$};
\node[align=left] (punto) at (5,-0) {$1$};
\draw (-6.25,-1.2) circle (0.3cm);
\draw (-5,-1.2) circle (0.3cm);
\draw (-3.75,-1.2) circle (0.3cm);
\draw (-2.5,-1.2) circle (0.3cm);
\draw (-1.25,-1.2) circle (0.3cm);
\draw (0,-1.2) circle (0.3cm);
\draw (1.25,-1.2) circle (0.3cm);
\draw (2.5,-1.2) circle (0.3cm);
\draw (3.75,-1.2) circle (0.3cm);
\draw (5,-1.2) circle (0.3cm);
\draw[thick,-] (-5.95,-1.2) -- (-5.35,-1.2);
\draw[thick,-] (-4.65,-1.2) -- (-4.1,-1.2);
\draw[thick,-] (-3.4,-1.2) -- (-2.85,-1.2);
\draw[thick,-] (-2.15,-1.2) -- (-1.6,-1.2);
\draw[thick,-] (-0.9,-1.2) -- (-0.35,-1.2);
\draw[thick,-] (0.35,-1.2) -- (0.9,-1.2);
\draw[thick,-] (1.6,-1.2) -- (2.15,-1.2);
\draw[thick,-] (2.85,-1.2) -- (3.4,-1.2);
\draw[thick,-] (4.1,-1.2) -- (4.65,-1.2);
\node[align=left] (punto) at (-8.5,-1.2) {${\mu}_2=\frac 19,$};
\node[align=left] (punto) at (-7.25,-1.2) {${\bf f}_2:\frac 1{18}$};
\node[align=left] (punto) at (-6.25,-1.2) {$1$};
\node[align=left] (punto) at (-5,-1.2) {$1$};
\node[align=left] (punto) at (-3.75,-1.2) {$1$};
\node[align=left] (punto) at (-2.5,-1.2) {$1$};
\node[align=left] (punto) at (-1.25,-1.2) {$1$};
\node[align=left] (punto) at (0,-1.2) {$-1$};
\node[align=left] (punto) at (1.25,-1.2) {$-1$};
\node[align=left] (punto) at (2.5,-1.2) {$-1$};
\node[align=left] (punto) at (3.75,-1.2) {$-1$};
\node[align=left] (punto) at (5,-1.2) {$-1$};
\draw (-6.25,-2.4) circle (0.3cm);
\draw (-5,-2.4) circle (0.3cm);
\draw (-3.75,-2.4) circle (0.3cm);
\draw (-2.5,-2.4) circle (0.3cm);
\draw (-1.25,-2.4) circle (0.3cm);
\draw (0,-2.4) circle (0.3cm);
\draw (1.25,-2.4) circle (0.3cm);
\draw (2.5,-2.4) circle (0.3cm);
\draw (3.75,-2.4) circle (0.3cm);
\draw (5,-2.4) circle (0.3cm);
\draw[thick,-] (-5.95,-2.4) -- (-5.35,-2.4);
\draw[thick,-] (-4.65,-2.4) -- (-4.1,-2.4);
\draw[thick,-] (-3.4,-2.4) -- (-2.85,-2.4);
\draw[thick,-] (-2.15,-2.4) -- (-1.6,-2.4);
\draw[thick,-] (-0.9,-2.4) -- (-0.35,-2.4);
\draw[thick,-] (0.35,-2.4) -- (0.9,-2.4);
\draw[thick,-] (1.6,-2.4) -- (2.15,-2.4);
\draw[thick,-] (2.85,-2.4) -- (3.4,-2.4);
\draw[thick,-] (4.1,-2.4) -- (4.65,-2.4);
\node[align=left] (punto) at (-8.5,-2.4) {${\mu}_3=\frac 17,$};
\node[align=left] (punto) at (-7.25,-2.4) {${\bf f}_3:\frac1{14}$};
\node[align=left] (punto) at (-6.25,-2.4) {$1$};
\node[align=left] (punto) at (-5,-2.4) {$1$};
\node[align=left] (punto) at (-3.75,-2.4) {$1$};
\node[align=left] (punto) at (-2.5,-2.4) {$1$};
\node[align=left] (punto) at (-1.25,-2.4) {$0$};
\node[align=left] (punto) at (0,-2.4) {$0$};
\node[align=left] (punto) at (1.25,-2.4) {$-1$};
\node[align=left] (punto) at (2.5,-2.4) {$-1$};
\node[align=left] (punto) at (3.75,-2.4) {$-1$};
\node[align=left] (punto) at (5,-2.4) {$-1$};
\draw (-6.25,-3.6) circle (0.3cm);
\draw (-5,-3.6) circle (0.3cm);
\draw (-3.75,-3.6) circle (0.3cm);
\draw (-2.5,-3.6) circle (0.3cm);
\draw (-1.25,-3.6) circle (0.3cm);
\draw (0,-3.6) circle (0.3cm);
\draw (1.25,-3.6) circle (0.3cm);
\draw (2.5,-3.6) circle (0.3cm);
\draw (3.75,-3.6) circle (0.3cm);
\draw (5,-3.6) circle (0.3cm);
\draw[thick,-] (-5.95,-3.6) -- (-5.35,-3.6);
\draw[thick,-] (-4.65,-3.6) -- (-4.1,-3.6);
\draw[thick,-] (-3.4,-3.6) -- (-2.85,-3.6);
\draw[thick,-] (-2.15,-3.6) -- (-1.6,-3.6);
\draw[thick,-] (-0.9,-3.6) -- (-0.35,-3.6);
\draw[thick,-] (0.35,-3.6) -- (0.9,-3.6);
\draw[thick,-] (1.6,-3.6) -- (2.15,-3.6);
\draw[thick,-] (2.85,-3.6) -- (3.4,-3.6);
\draw[thick,-] (4.1,-3.6) -- (4.65,-3.6);
\node[align=left] (punto) at (-8.5,-3.6) {${\mu}_4=\frac 15,$};
\node[align=left] (punto) at (-7.25,-3.6) {${\bf f}_4:\frac1{10}$};
\node[align=left] (punto) at (-6.25,-3.6) {$1$};
\node[align=left] (punto) at (-5,-3.6) {$1$};
\node[align=left] (punto) at (-3.75,-3.6) {$1$};
\node[align=left] (punto) at (-2.5,-3.6) {$0$};
\node[align=left] (punto) at (-1.25,-3.6) {$0$};
\node[align=left] (punto) at (0,-3.6) {$0$};
\node[align=left] (punto) at (1.25,-3.6) {$0$};
\node[align=left] (punto) at (2.5,-3.6) {$-1$};
\node[align=left] (punto) at (3.75,-3.6) {$-1$};
\node[align=left] (punto) at (5,-3.6) {$-1$};
\draw (-6.25,-4.8) circle (0.3cm);
\draw (-5,-4.8) circle (0.3cm);
\draw (-3.75,-4.8) circle (0.3cm);
\draw (-2.5,-4.8) circle (0.3cm);
\draw (-1.25,-4.8) circle (0.3cm);
\draw (0,-4.8) circle (0.3cm);
\draw (1.25,-4.8) circle (0.3cm);
\draw (2.5,-4.8) circle (0.3cm);
\draw (3.75,-4.8) circle (0.3cm);
\draw (5,-4.8) circle (0.3cm);
\draw[thick,-] (-5.95,-4.8) -- (-5.35,-4.8);
https://it.overleaf.com/project/5f047e63e80c230001a9e222\draw[thick,-] (-4.65,-4.8) -- (-4.1,-4.8);
\draw[thick,-] (-3.4,-4.8) -- (-2.85,-4.8);
\draw[thick,-] (-2.15,-4.8) -- (-1.6,-4.8);
\draw[thick,-] (-0.9,-4.8) -- (-0.35,-4.8);
\draw[thick,-] (0.35,-4.8) -- (0.9,-4.8);
\draw[thick,-] (1.6,-4.8) -- (2.15,-4.8);
\draw[thick,-] (2.85,-4.8) -- (3.4,-4.8);
\draw[thick,-] (4.1,-4.8) -- (4.65,-4.8);
\node[align=left] (punto) at (-8.5,-4.8) {${\mu}_5=\frac 14,$};
\node[align=left] (punto) at (-7.25,-4.8) {${\bf f}_5:\frac1{8}$};
\node[align=left] (punto) at (-6.25,-4.8) {$0$};
\node[align=left] (punto) at (-5,-4.8) {$0$};
\node[align=left] (punto) at (-3.75,-4.8) {$0$};
\node[align=left] (punto) at (-2.5,-4.8) {$1$};
\node[align=left] (punto) at (-1.25,-4.8) {$1$};
\node[align=left] (punto) at (0,-4.8) {$1$};
\node[align=left] (punto) at (1.25,-4.8) {$1$};
\node[align=left] (punto) at (2.5,-4.8) {$0$};
\node[align=left] (punto) at (3.75,-4.8) {$0$};
\node[align=left] (punto) at (5,-4.8) {$0$};
\draw (-6.25,-6) circle (0.3cm);
\draw (-5,-6) circle (0.3cm);
\draw (-3.75,-6) circle (0.3cm);
\draw (-2.5,-6) circle (0.3cm);
\draw (-1.25,-6) circle (0.3cm);
\draw (0,-6) circle (0.3cm);
\draw (1.25,-6) circle (0.3cm);
\draw (2.5,-6) circle (0.3cm);
\draw (3.75,-6) circle (0.3cm);
\draw (5,-6) circle (0.3cm);
\draw[thick,-] (-5.95,-6) -- (-5.35,-6);
\draw[thick,-] (-4.65,-6) -- (-4.1,-6);
\draw[thick,-] (-3.4,-6) -- (-2.85,-6);
\draw[thick,-] (-2.15,-6) -- (-1.6,-6);
\draw[thick,-] (-0.9,-6) -- (-0.35,-6);
\draw[thick,-] (0.35,-6) -- (0.9,-6);
\draw[thick,-] (1.6,-6) -- (2.15,-6);
\draw[thick,-] (2.85,-6) -- (3.4,-6);
\draw[thick,-] (4.1,-6) -- (4.65,-6);
\node[align=left] (punto) at (-8.5,-6) {${\mu}_6=\frac 13,$};
\node[align=left] (punto) at (-7.25,-6) {${\bf f}_6:\frac1{18}$};
\node[align=left] (punto) at (-6.25,-6) {$1$};
\node[align=left] (punto) at (-5,-6) {$1$};
\node[align=left] (punto) at (-3.75,-6) {$-1$};
\node[align=left] (punto) at (-2.5,-6) {$-1$};
\node[align=left] (punto) at (-1.25,-6) {$-1$};
\node[align=left] (punto) at (0,-6) {$1$};
\node[align=left] (punto) at (1.25,-6) {$1$};
\node[align=left] (punto) at (2.5,-6) {$-1$};
\node[align=left] (punto) at (3.75,-6) {$-1$};
\node[align=left] (punto) at (5,-6) {$-1$};
\draw (-6.25,-7.2) circle (0.3cm);
\draw (-5,-7.2) circle (0.3cm);
\draw (-3.75,-7.2) circle (0.3cm);
\draw (-2.5,-7.2) circle (0.3cm);
\draw (-1.25,-7.2) circle (0.3cm);
\draw (0,-7.2) circle (0.3cm);
\draw (1.25,-7.2) circle (0.3cm);
\draw (2.5,-7.2) circle (0.3cm);
\draw (3.75,-7.2) circle (0.3cm);
\draw (5,-7.2) circle (0.3cm);
\draw[thick,-] (-5.95,-7.2) -- (-5.35,-7.2);
\draw[thick,-] (-4.65,-7.2) -- (-4.1,-7.2);
\draw[thick,-] (-3.4,-7.2) -- (-2.85,-7.2);
\draw[thick,-] (-2.15,-7.2) -- (-1.6,-7.2);
\draw[thick,-] (-0.9,-7.2) -- (-0.35,-7.2);
\draw[thick,-] (0.35,-7.2) -- (0.9,-7.2);
\draw[thick,-] (1.6,-7.2) -- (2.15,-7.2);
\draw[thick,-] (2.85,-7.2) -- (3.4,-7.2);
\draw[thick,-] (4.1,-7.2) -- (4.65,-7.2);
\node[align=left] (punto) at (-8.5,-7.2) {${\mu}_7=\frac 12,$};
\node[align=left] (punto) at (-7.25,-7.2) {${\bf f}_7:\frac1{12}$};
\node[align=left] (punto) at (-6.25,-7.2) {$0$};
\node[align=left] (punto) at (-5,-7.2) {$1$};
\node[align=left] (punto) at (-3.75,-7.2) {$1$};
\node[align=left] (punto) at (-2.5,-7.2) {$0$};
\node[align=left] (punto) at (-1.25,-7.2) {$-1$};
\node[align=left] (punto) at (0,-7.2) {$-1$};
\node[align=left] (punto) at (1.25,-7.2) {$0$};
\node[align=left] (punto) at (2.5,-7.2) {$1$};
\node[align=left] (punto) at (3.75,-7.2) {$1$};
\node[align=left] (punto) at (5,-7.2) {$0$};
\draw (-6.25,-8.4) circle (0.3cm);
\draw (-5,-8.4) circle (0.3cm);
\draw (-3.75,-8.4) circle (0.3cm);
\draw (-2.5,-8.4) circle (0.3cm);
\draw (-1.25,-8.4) circle (0.3cm);
\draw (0,-8.4) circle (0.3cm);
\draw (1.25,-8.4) circle (0.3cm);
\draw (2.5,-8.4) circle (0.3cm);
\draw (3.75,-8.4) circle (0.3cm);
\draw (5,-8.4) circle (0.3cm);
\draw[thick,-] (-5.95,-8.4) -- (-5.35,-8.4);
\draw[thick,-] (-4.65,-8.4) -- (-4.1,-8.4);
\draw[thick,-] (-3.4,-8.4) -- (-2.85,-8.4);
\draw[thick,-] (-2.15,-8.4) -- (-1.6,-8.4);
\draw[thick,-] (-0.9,-8.4) -- (-0.35,-8.4);
\draw[thick,-] (0.35,-8.4) -- (0.9,-8.4);
\draw[thick,-] (1.6,-8.4) -- (2.15,-8.4);
\draw[thick,-] (2.85,-8.4) -- (3.4,-8.4);
\draw[thick,-] (4.1,-8.4) -- (4.65,-8.4);
\node[align=left] (punto) at (-8.5,-8.4) {${\mu}_8=1,$};
\node[align=left] (punto) at (-7.25,-8.4) {${\bf f}_8:\frac1{18}$};
\node[align=left] (punto) at (-6.25,-8.4) {$1$};
\node[align=left] (punto) at (-5,-8.4) {$-1$};
\node[align=left] (punto) at (-3.75,-8.4) {$1$};
\node[align=left] (punto) at (-2.5,-8.4) {$-1$};
\node[align=left] (punto) at (-1.25,-8.4) {$1$};
\node[align=left] (punto) at (0,-8.4) {$-1$};
\node[align=left] (punto) at (1.25,-8.4) {$1$};
\node[align=left] (punto) at (2.5,-8.4) {$-1$};
\node[align=left] (punto) at (3.75,-8.4) {$1$};
\node[align=left] (punto) at (5,-8.4) {$-1$};
\end{tikzpicture}
\end{center}
The second eigenvector ${\bf f}_2$ is orthogonal to ${\bf f}_1$, according with the results on the second Cheeger constant. It is easily seen that ${\bf f}_5$ achieves the lower eigenvalue among the eigenvectors pseudo-orthogonal to ${\bf f}_1$ and ${\bf f}_2$. Furthermore, we have $h_3(P_{10})=\mu_5(P_{10})=\frac 14$. 
 \end{exa}
 \begin{exa}
The cycle graph $C_{10}$ is a graph for which every vertices have degree $2$. The spectrum is $\sigma(C_{10})=\left\{0,\frac 15, \frac 14, \frac 13, \frac 12,1\right\}$.
\begin{center}
\begin{tikzpicture}
\node[align=left] (punto) at (-8,2) {${\mu}_1=0$};
\node[align=left] (punto) at (-8,1) {${\bf f}_1:\frac1{20}$};
\node[align=left] (punto) at (-5.5,2) {${\mu}_2=\frac 15$};
\node[align=left] (punto) at (-5.5,1) {${\bf f}_2:\frac 1{20}$};
\node[align=left] (punto) at (-3,2) {${\mu}_3=\frac 14$};
\node[align=left] (punto) at (-3.25,1) {${\bf f}_3:\frac1{16}$};
\node[align=left] (punto) at (-0.5,2) {${\mu}_4=\frac 13$};
\node[align=left] (punto) at (-0.5,1) {${\bf f}_4:\frac1{12}$};
\node[align=left] (punto) at (2,2) {${\mu}_5=\frac 12$};
\node[align=left] (punto) at (2,1) {${\bf f}_5:\frac1{16}$};
\node[align=left] (punto) at (4.5,2) {${\mu}_6=1$};
\node[align=left] (punto) at (4.5,1) {${\bf f}_6:\frac1{20}$};
\draw (-8.75,0) circle (0.3cm);
\draw (-7.5,0) circle (0.3cm);
\draw (-6.25,0) circle (0.3cm);
\draw (-5,0) circle (0.3cm);
\draw (-3.75,0) circle (0.3cm);
\draw (-2.5,0) circle (0.3cm);
\draw (-1.25,0) circle (0.3cm);
\draw (0,0) circle (0.3cm);
\draw (1.25,0) circle (0.3cm);
\draw (2.5,0) circle (0.3cm);
\draw (3.75,0) circle (0.3cm);
\draw (5,0) circle (0.3cm);
\draw[thick,-] (-8.45,0) -- (-7.85,0);
\draw[thick,-] (-5.95,0) -- (-5.35,0);
\draw[thick,-] (-3.4,0) -- (-2.85,0);
\draw[thick,-] (-0.9,0) -- (-0.35,0);
\draw[thick,-] (1.6,0) -- (2.15,0);
\draw[thick,-] (4.1,0) -- (4.65,0);
\draw[thick,-] (-8.75,-0.35) -- (-8.75,-0.85);
\draw[thick,-] (-7.5,-0.35) -- (-7.5,-0.85);
\draw[thick,-] (-6.25,-0.35) -- (-6.25,-0.85);
\draw[thick,-] (-5,-0.35) -- (-5,-0.85);
\draw[thick,-] (-3.75,-0.35) -- (-3.75,-0.85);
\draw[thick,-] (-2.5,-0.35) -- (-2.5,-0.85);
\draw[thick,-] (-1.25,-0.35) -- (-1.25,-0.85);
\draw[thick,-] (0,-0.35) -- (0,-0.85);
\draw[thick,-] (1.25,-0.35) -- (1.25,-0.85);
\draw[thick,-] (2.5,-0.35) -- (2.5,-0.85);
\draw[thick,-] (3.75,-0.35) -- (3.75,-0.85);
\draw[thick,-] (5,-0.35) -- (5,-0.85);
\node[align=left] (punto) at (-8.5,0.5) {$1$};
\node[align=left] (punto) at (-7.25,0.5) {$2$};
\node[align=left] (punto) at (-8.5,-0.75) {$10$};
\node[align=left] (punto) at (-7.25,-0.75) {$3$};
\node[align=left] (punto) at (-8.5,-1.9) {$9$};
\node[align=left] (punto) at (-7.25,-1.9) {$4$};
\node[align=left] (punto) at (-8.5,-3.15) {$8$};
\node[align=left] (punto) at (-7.25,-3.15) {$5$};
\node[align=left] (punto) at (-8.5,-4.35) {$7$};
\node[align=left] (punto) at (-7.25,-4.35) {$6$};
\node[align=left] (punto) at (-8.75,-0) {$1$};
\node[align=left] (punto) at (-7.5,-0) {$1$};
\node[align=left] (punto) at (-6.25,-0) {$1$};
\node[align=left] (punto) at (-5,-0) {$-1$};
\node[align=left] (punto) at (-3.75,-0) {$0$};
\node[align=left] (punto) at (-2.5,0) {$0$};
\node[align=left] (punto) at (-1.25,0) {$0$};
\node[align=left] (punto) at (0,0) {$0$};
\node[align=left] (punto) at (1.25,0) {$0$};
\node[align=left] (punto) at (2.5,-0) {$0$};
\node[align=left] (punto) at (3.75,-0) {$1$};
\node[align=left] (punto) at (5,-0) {$-1$};
\draw (-8.75,-1.2) circle (0.3cm);
\draw (-7.5,-1.2) circle (0.3cm);
\draw (-6.25,-1.2) circle (0.3cm);
\draw (-5,-1.2) circle (0.3cm);
\draw (-3.75,-1.2) circle (0.3cm);
\draw (-2.5,-1.2) circle (0.3cm);
\draw (-1.25,-1.2) circle (0.3cm);
\draw (0,-1.2) circle (0.3cm);
\draw (1.25,-1.2) circle (0.3cm);
\draw (2.5,-1.2) circle (0.3cm);
\draw (3.75,-1.2) circle (0.3cm);
\draw (5,-1.2) circle (0.3cm);
\draw[thick,-] (-8.75,-1.55) -- (-8.75,-2.05);
\draw[thick,-] (-7.5,-1.55) -- (-7.5,-2.05);
\draw[thick,-] (-6.25,-1.55) -- (-6.25,-2.05);
\draw[thick,-] (-5,-1.55) -- (-5,-2.05);
\draw[thick,-] (-3.75,-1.55) -- (-3.75,-2.05);
\draw[thick,-] (-2.5,-1.55) -- (-2.5,-2.05);
\draw[thick,-] (-1.25,-1.55) -- (-1.25,-2.05);
\draw[thick,-] (0,-1.55) -- (0,-2.05);
\draw[thick,-] (1.25,-1.55) -- (1.25,-2.05);
\draw[thick,-] (2.5,-1.55) -- (2.5,-2.05);
\draw[thick,-] (3.75,-1.55) -- (3.75,-2.05);
\draw[thick,-] (5,-1.55) -- (5,-2.05);
\node[align=left] (punto) at (-8.75,-1.2) {$1$};
\node[align=left] (punto) at (-7.5,-1.2) {$1$};
\node[align=left] (punto) at (-6.25,-1.2) {$1$};
\node[align=left] (punto) at (-5,-1.2) {$-1$};
\node[align=left] (punto) at (-3.75,-1.2) {$1$};
\node[align=left] (punto) at (-2.5,-1.2) {$-1$};
\node[align=left] (punto) at (-1.25,-1.2) {$0$};
\node[align=left] (punto) at (0,-1.2) {$0$};
\node[align=left] (punto) at (1.25,-1.2) {$-1$};
\node[align=left] (punto) at (2.5,-1.2) {$1$};
\node[align=left] (punto) at (3.75,-1.2) {$-1$};
\node[align=left] (punto) at (5,-1.2) {$1$};
\draw (-8.75,-2.4) circle (0.3cm);
\draw (-7.5,-2.4) circle (0.3cm);
\draw (-6.25,-2.4) circle (0.3cm);
\draw (-5,-2.4) circle (0.3cm);
\draw (-3.75,-2.4) circle (0.3cm);
\draw (-2.5,-2.4) circle (0.3cm);
\draw (-1.25,-2.4) circle (0.3cm);
\draw (0,-2.4) circle (0.3cm);
\draw (1.25,-2.4) circle (0.3cm);
\draw (2.5,-2.4) circle (0.3cm);
\draw (3.75,-2.4) circle (0.3cm);
\draw (5,-2.4) circle (0.3cm);
\node[align=left] (punto) at (-8.75,-2.4) {$1$};
\node[align=left] (punto) at (-7.5,-2.4) {$1$};
\node[align=left] (punto) at (-6.25,-2.4) {$1$};
\node[align=left] (punto) at (-5,-2.4) {$-1$};
\node[align=left] (punto) at (-3.75,-2.4) {$1$};
\node[align=left] (punto) at (-2.5,-2.4) {$-1$};
\node[align=left] (punto) at (-1.25,-2.4) {$1$};
\node[align=left] (punto) at (0,-2.4) {$-1$};
\node[align=left] (punto) at (1.25,-2.4) {$-1$};
\node[align=left] (punto) at (2.5,-2.4) {$1$};
\node[align=left] (punto) at (3.75,-2.4) {$1$};
\node[align=left] (punto) at (5,-2.4) {$-1$};
\draw (-8.75,-3.6) circle (0.3cm);
\draw (-7.5,-3.6) circle (0.3cm);
\draw (-6.25,-3.6) circle (0.3cm);
\draw (-5,-3.6) circle (0.3cm);
\draw (-3.75,-3.6) circle (0.3cm);
\draw (-2.5,-3.6) circle (0.3cm);
\draw (-1.25,-3.6) circle (0.3cm);
\draw (0,-3.6) circle (0.3cm);
\draw (1.25,-3.6) circle (0.3cm);
\draw (2.5,-3.6) circle (0.3cm);
\draw (3.75,-3.6) circle (0.3cm);
\draw (5,-3.6) circle (0.3cm);
\draw[thick,-] (-8.75,-2.75) -- (-8.75,-3.25);
\draw[thick,-] (-7.5,-2.75) -- (-7.5,-3.25);
\draw[thick,-] (-6.25,-2.75) -- (-6.25,-3.25);
\draw[thick,-] (-5,-2.75) -- (-5,-3.25);
\draw[thick,-] (-3.75,-2.75) -- (-3.75,-3.25);
\draw[thick,-] (-2.5,-2.75) -- (-2.5,-3.25);
\draw[thick,-] (-1.25,-2.75) -- (-1.25,-3.25);
\draw[thick,-] (0,-2.75) -- (0,-3.25);
\draw[thick,-] (1.25,-2.75) -- (1.25,-3.25);
\draw[thick,-] (2.5,-2.75) -- (2.5,-3.25);
\draw[thick,-] (3.75,-2.75) -- (3.75,-3.25);
\draw[thick,-] (5,-2.75) -- (5,-3.25);
\node[align=left] (punto) at (-8.75,-3.6) {$1$};
\node[align=left] (punto) at (-7.5,-3.6) {$1$};
\node[align=left] (punto) at (-6.25,-3.6) {$1$};
\node[align=left] (punto) at (-5,-3.6) {$-1$};
\node[align=left] (punto) at (-3.75,-3.6) {$1$};
\node[align=left] (punto) at (-2.5,-3.6) {$-1$};
\node[align=left] (punto) at (-1.25,-3.6) {$1$};
\node[align=left] (punto) at (0,-3.6) {$-1$};
\node[align=left] (punto) at (1.25,-3.6) {$1$};
\node[align=left] (punto) at (2.5,-3.6) {$-1$};
\node[align=left] (punto) at (3.75,-3.6) {$-1$};
\node[align=left] (punto) at (5,-3.6) {$1$};
\draw (-8.75,-4.8) circle (0.3cm);
\draw (-7.5,-4.8) circle (0.3cm);
\draw (-6.25,-4.8) circle (0.3cm);
\draw (-5,-4.8) circle (0.3cm);
\draw (-3.75,-4.8) circle (0.3cm);
\draw (-2.5,-4.8) circle (0.3cm);
\draw (-1.25,-4.8) circle (0.3cm);
\draw (0,-4.8) circle (0.3cm);
\draw (1.25,-4.8) circle (0.3cm);
\draw (2.5,-4.8) circle (0.3cm);
\draw (3.75,-4.8) circle (0.3cm);
\draw (5,-4.8) circle (0.3cm);
\node[align=left] (punto) at (-8.75,-4.8) {$1$};
\node[align=left] (punto) at (-7.5,-4.8) {$1$};
\draw[thick,-] (-5.95,-4.8) -- (-5.35,-4.8);
\draw[thick,-] (-3.4,-4.8) -- (-2.85,-4.8);
\draw[thick,-] (-0.9,-4.8) -- (-0.35,-4.8);
\draw[thick,-] (1.6,-4.8) -- (2.15,-4.8);
\draw[thick,-] (4.1,-4.8) -- (4.65,-4.8);
\draw[thick,-] (-8.4,-4.8) -- (-7.85,-4.8);
\draw[thick,-] (-8.75,-3.95) -- (-8.75,-4.45);
\draw[thick,-] (-7.5,-3.95) -- (-7.5,-4.45);
\draw[thick,-] (-6.25,-3.95) -- (-6.25,-4.45);
\draw[thick,-] (-5,-3.95) -- (-5,-4.45);
\draw[thick,-] (-3.75,-3.95) -- (-3.75,-4.45);
\draw[thick,-] (-2.5,-3.95) -- (-2.5,-4.45);
\draw[thick,-] (-1.25,-3.95) -- (-1.25,-4.45);
\draw[thick,-] (0,-3.95) -- (0,-4.45);
\draw[thick,-] (1.25,-3.95) -- (1.25,-4.45);
\draw[thick,-] (2.5,-3.95) -- (2.5,-4.45);
\draw[thick,-] (3.75,-3.95) -- (3.75,-4.45);
\draw[thick,-] (5,-3.95) -- (5,-4.45);
\node[align=left] (punto) at (-6.25,-4.8) {$1$};
\node[align=left] (punto) at (-5,-4.8) {$-1$};
\node[align=left] (punto) at (-3.75,-4.8) {$1$};
\node[align=left] (punto) at (-2.5,-4.8) {$-1$};
\node[align=left] (punto) at (-1.25,-4.8) {$1$};
\node[align=left] (punto) at (0,-4.8) {$-1$};
\node[align=left] (punto) at (1.25,-4.8) {$1$};
\node[align=left] (punto) at (2.5,-4.8) {$-1$};
\node[align=left] (punto) at (3.75,-4.8) {$1$};
\node[align=left] (punto) at (5,-4.8) {$-1$};
\end{tikzpicture}
\end{center} 
The second eigenvector ${\bf f}_2$ is orthogonal to ${\bf f}_1$, according with the results on the second Cheeger constant. It is easily seen and that ${\bf f}_4$ (to be more precise, only the postive part) achieves the lower eigenvalue among the eigenvectors pseudo orthogonal to ${\bf f}_1$ and ${\bf f}_2$. Furthermore, we have $h_3(P_{10})=\mu_4(P_{10})=\frac 14$. 
 \end{exa}

\section{The Pseudo-orthogonality for characterizing the Cheeger Constants }\label{aSCCC}
Throughout this paper we assume that it is always possible to consider three disjoint non-empty sets of $V$. In this way, it is always possible to define the third Cheeger constant. Hence, we focus on giving some characterizations of
non-trivial 1-Laplacian eigenvalues in the form of continuous optimizations.

\subsection{The asymptotic behaviour of the graph p-Laplacian eigenvalue in the continuous case}
To motivate our treatment on the graph $1$-Laplacian, we remark that it is deeply related to the Cheeger problem, also in the continuous case. Indeed, let $\Omega$ be a bounded domain of $\R^{n}$, $n\ge 2$, then, for any $1<p<+\infty$, the $p$-Laplacian operator is defined as
\begin{equation*}
 \Delta_{p}u:=\dive \left(|\nabla u|^{p-2}\nabla u\right),\quad u\in W^{1,p}(\Omega).
\end{equation*} 
Let us consider the following Dirichlet eigenvalue problem:
\begin{equation}
\label{introauti}
\left\{
\begin{array}{ll}
-\Delta_{p}u=\lambda^{(p)}(\Omega) |u|^{p-2}u  &\text{in}\ \Omega, \\
u=0 &\text{on } \partial\Omega.
\end{array}
\right.
\end{equation}
Cheeger \cite{Che} proved the \lq\lq Cheeger inequality\rq\rq:
\begin{equation*}
\lambda^{(p)}(\Omega)\geq\left(\frac{\rho_{1}(\Omega)}{p}\right)^p,
\end{equation*}
where $\rho_{1}(\Omega)$ is the first k-way Cheeger constant, defined as
\begin{equation*}
\rho_{1}(\Omega):=\inf_{\substack{E\subset\Omega}}\left\{\frac{P(E)}{|E|}\right\}.
\end{equation*}
with $P$ denoting the perimeter of $E$ in $\R^n$. Afterwards, Kawohl and Fridman \cite{KF} studied the asymptotic behavior of the first eigenvalue of \eqref{introauti}, as $p\to1$:
\begin{equation}
\label{firsttends}
\lim_{p\to 1}\lambda_1^{(p)}(\Omega)^{\frac1p}=\rho_{1}(\Omega).
\end{equation}
In \cite{LS}, the authors show the asymptotic convergence to the $1$-Laplace eigenvalues:
\begin{equation*}
\lim_{p\to 1}\lambda_k^{(p)}(\Omega)^{\frac1p}=\lambda_k^{(1)}(\Omega)\quad\forall k\in\N.
\end{equation*}
Then, similarly to \eqref{k-way-cheeger-constant}, in the continuous case, the higher $k$-way Cheeger constants are defined as
\begin{equation*}
\rho_{k}(\Omega):=\inf \left\{ \max_{i=1,...,k}\frac{P(E_{i})}{|E_{i}|}\ :\ E_{i}\subset \Omega,\ |E_{i}|>0\ \forall i, \ E_{i}\cap E_{j}=\emptyset\ \forall i\neq j\right\}.
\end{equation*}
Subsequently, asymptotic results as in \eqref{firsttends} have been generalized to higher eigenvalues. 

Particularly, in \cite{Pa}, it has been proven the Cheeger inequality for the second eigenvalue $\lambda_{2}^{(p)}(\Omega)$ of \eqref{introauti}:
\begin{equation*}
\lambda_{2}^{(p)}(\Omega)\geq\left(\frac{\rho_{2}(\Omega)}{p}\right)^p
\end{equation*} 
and the limit property 
\begin{equation}\label{asihig}
\limsup_{p\to 1}\lambda_2^{(p)}(\Omega)^{\frac1p}=\rho_{2}(\Omega) \quad\text{and}\quad\limsup_{p\to 1}\lambda_k^{(p)}(\Omega)^{\frac1p}\leq \rho_{k}(\Omega).
\end{equation}
The reason of the discrepancy between the two relations in \eqref{asihig} relies on the fact that every second eigenfunctions  has exactly two nodal domains but, on the other hand, a $k$-th eigenfunction generally have not $k$ nodal domains.

Regarding the higher Cheeger constant, it has been proven that the two quantities
\begin{equation*}
\begin{split}
&\Lambda_k^{(p)}(\Omega):=\inf \left\{ \sum_{i=1,...,k}\lambda_1^{(p)}(E_{i})\ :\ E_{i}\subset \Omega,\ |E_{i}|>0\ \forall i, \ E_{i}\cap E_{j}=\emptyset\ \forall i\neq j\right\},\\
&\mathcal L_k(p,\Omega):=\inf \left\{ \max_{i=1,...,k}\lambda_1^{(p)}(E_{i})\ :\ E_{i}\subset \Omega,\ |E_{i}|>0\ \forall i, \ E_{i}\cap E_{j}=\emptyset\ \forall i\neq j\right\},
\end{split}
\end{equation*}
have the following limit behaviours (see \cite{Ca} and \cite{BP} for further details):
\begin{equation*}
\begin{split}
&\lim_{p\to 1}\Lambda_k^{(p)}=
\inf \left\{ \sum_{i=1,...,k}\rho_1(p,E_{i})\ :\ E_{i}\subset \Omega,\ |E_{i}|>0\ \forall i, \ E_{i}\cap E_{j}=\emptyset\ \forall i\neq j\right\},\\
&\lim_{p\to 1}\mathcal L_{k}(p,\Omega)=\rho_k(\Omega).
\end{split}
\end{equation*}

The relationship between the $1$-Laplacian eigenvalues and the Cheeger constants has been investigated also in the discrete case. Indeed, the graph $1$-Laplacian is the limiting operator, as $p\to1$, of the graph $p$-Laplacian:
\begin{equation}
\label{graph_p_lap}
(\Delta_p {\bf f})_i:=\sum_{\substack{i,j\in V\\ j\sim i}} |f_i-f_j|^{p-2}(f_i-f_j)
\end{equation}
and it has been proven that the second eigenvalue of the graph $p$-Laplacian  approximate the second Cheeger constant arbitrarily well (see \cite{A,BHa}). 

Regarding the approximation of the higher Cheeger constants with the higher eigenvalues of the $p$-Laplacian \eqref{graph_p_lap} as $p\to 1$, the following Cheeger inequality (see \cite{TH}) holds for any $k\in\N$:
\begin{equation}
\label{asy_eig}
\left(\dfrac 2 {\max_i d_i}\right)^{p-1}\dfrac{\rho_{j}(G)}{p^p}^p\leq\mu_k^{(p)}(G)\leq 2^{p-1}\rho_{k}(G),
\end{equation}
where $j=2,..,k$ is the number of nodal domains of the $k$-th eigenfunction.

Since the discrete nodal domain Theorem \cite{DLS,TH} states that the $k$-th eigenfunctions have at most $k$ nodal domains, than the inequality \eqref{asy_eig} gives better estimates when considering eigenvalues admitting eigenfunctions with exactly $k$ nodal domains.

\subsection{The role of the orthogonality}
In \cite[p.6]{Chu} the following characterization of the graph ($2$-)Laplacian eigenvalues is given:
\begin{equation}
\label{chucar}
\mu_k^{(2)}(G)=\min_{{\mathbf f}\perp C_{k-1}}\dfrac{\sum_{\substack{i,j\in V\\ j\sim i}}|f_i-f_j|^2}{\sum_{i \in V}d_i|f_i|^2}=\min_{{\mathbf f}\neq 0}\max_{{\mathbf v} \in C_{k-1}}\dfrac{\sum_{\substack{i,j\in V\\ j\sim i}}|f_i-f_j|^2}{\sum_{i \in V}d_i|f_i-v_i|^2},
\end{equation}
where $C_k$ is the subspace spanned by eigenfunctions achieving $\mu^{(2)}_j(G)$, for $1 \leq j \leq k$.

The fact that the $1$-Laplacian eigenvalues are asymptotically the Cheeger constants \eqref{asy_eig} and the characterizations in \eqref{chucar} motivate us to look for similar characterizations for the second, the third (and higher) Cheeger constants.

The following example motivates the use of the pseudo-orthogonality $\perp_p$. Indeed it is possible to find two different eigenvectors, associated to two different eigenvalues, such that they are pseudo-orthogonal but not orthogonal, in the sense that their scalar product is not zero. To this aim, let us consider the following. 
\begin{exa}
Let us consider $G=(V,E)$, where $V=\{1,2,3,4\}$ and $E=\{e_1=(1,2); e_2=(2,3); e_3=(3,4)\}$. The eigenvectors $\hat{\bf f}_2=\frac{1}3(1,1,0,0)$ and $\hat{\bf f}_3=(1,0,0,0)$, represented below, are, respectively, associated to the eigenvalues $\mu_2=\frac 13$ and $\mu_3=1$ and they are not orthogonal, since $\langle \hat{\bf f}_2,\hat{\bf f}_3\rangle=\frac 13$.
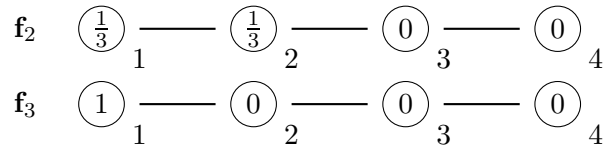
 \begin{center}
\begin{tikzpicture}
\draw[thick,-] (-1.5,0) -- (-0.5,0);
\draw[thick,-] (0.5,0) -- (1.5,0);
\draw[thick,-] (2.5,0) -- (3.5,0);
\draw (-2,0) circle (0.3cm);
\draw (0,0) circle (0.3cm);
\draw (2,0) circle (0.3cm);
\draw (4,0) circle (0.3cm);
\node[align=left] (punto) at (-1.5,-0.4) {$1$};
\node[align=left] (punto) at (0.5,-0.4) {$2$};
\node[align=left] (punto) at (2.5,-0.4) {$3$};
\node[align=left] (punto) at (4.5,-0.4) {$4$};
\node[align=left] (punto) at (-2,-0) {$\frac 13$};
\node[align=left] (punto) at (0,-0) {$\frac 13$};
\node[align=left] (punto) at (2,-0) {$0$};
\node[align=left] (punto) at (4,-0) {$0$};
\node[align=left] (punto) at (-3,-0) {${\bf f}_2$};
\end{tikzpicture}
\begin{tikzpicture}
\draw[thick,-] (-1.5,0) -- (-0.5,0);
\draw[thick,-] (0.5,0) -- (1.5,0);
\draw[thick,-] (2.5,0) -- (3.5,0);
\draw (-2,0) circle (0.3cm);
\draw (0,0) circle (0.3cm);
\draw (2,0) circle (0.3cm);
\draw (4,0) circle (0.3cm);
\node[align=left] (punto) at (-1.5,-0.4) {$1$};
\node[align=left] (punto) at (0.5,-0.4) {$2$};
\node[align=left] (punto) at (2.5,-0.4) {$3$};
\node[align=left] (punto) at (4.5,-0.4) {$4$};
\node[align=left] (punto) at (-2,-0) {$1$};
\node[align=left] (punto) at (0,-0) {$0$};
\node[align=left] (punto) at (2,-0) {$0$};
\node[align=left] (punto) at (4,-0) {$0$};
\node[align=left] (punto) at (-3,-0) {${\bf f}_3$};
\end{tikzpicture}
\end{center}  
\end{exa}

\subsection{The graph Cheeger constants in form of continuous optimizations}
There were developed many methods and techniques to clusterize a graph (refer to \cite{vL} for an overview) but, up to our knowledge, it is yet difficult to determine the optimal number of clusters in a data set, since it depends on the method used for measuring similarities and the parameter used for partitioning.
Since the $2$-clustering has been deeply studied (see e.g \cite{BHa,BHb} and reference therein), we focus on the $3$-clustering, that is the division of the nodes into three groups.
This is the motivation for which we focus on the third Cheeger constant. 

In this Section, for any $A, B\subseteq V$, we denote
\[
E(A, B):=\{(i,j)\in E\  | \ \text{either}\  i\in A,\ j\in B\ \text {or}\ j\in A,\ i\in B\}.\\
\]
In \cite[Th. 2.6]{Chu} and \cite[Lem. 5.14]{Cha} is given the following characterization of the second Cheeger constant. We improve the proof and generalize this characterization for the third Cheeger constant.
\begin{prop}
\label{espl_h2h3}
Let $G=(V,E)$ be a graph. Then there exist two vectors ${\bf y}_2$ and ${\bf y}_3$ such that
\begin{equation*}
\begin{split}
h_2(G)&=\max_{c\in\R} \dfrac{\sum_{\substack{i,j\in V\\ j\sim i}}|(y_2)_i-(y_2)_j|}{\sum_{i\in V}d_i|(y_2)_i-c|},\\
h_3(G)&\leq\max_{c_1, c_2 \in\R} \dfrac{\sum_{\substack{i,j\in V\\ j\sim i}}{|(y_3)_i-(y_3)_j-c_2((y_2)_i-(y_2)_j)|}}{\sum_{i\in V}d_i|(y_3)_i-c_1-c_2(y_2)_i|}.
\end{split}
\end{equation*}
\end{prop}
\begin{proof}
By definition, there exists a set $A\subseteq V$ such that
\begin{equation*}
h_2(G)=\dfrac {|\partial A|} {{\vol(A)}}\quad \text{with}\ \vol(A)\leq \vol(A^c),
\end{equation*}
where we have denoted $A^c=V\setminus A$. We verify that ${\bf y}_2:={\bf 1}_A$. Indeed,
\begin{equation*}
\begin{split}
\max_{c\in\R}\dfrac{\sum_{\substack{i,j\in V\\ j\sim i}}|(y_2)_i-(y_2)_j|}{ \sum_{i\in V}d_i|(y_2)_i-c|}
&=\max_{c\in\R} \dfrac{\sum_{\substack{i,j\in V\\ j\sim i}}|(  1_A)_i- (  1_A)_j|}{\sum_{i\in V}d_i|(  1_A)_i-c|}\\
&=\dfrac{\sum_{\substack{i,j\in V\\ j\sim i}}|(  1_A)_i-(  1_A)_j|}{\min_{0\leq c\le 1} \sum_{i\in V}d_i|(  1_A)_i-c|}\\
&= \dfrac{|\partial A|}{\min_{0\leq c\leq 1}(1-c)\vol(A)+ c\vol(A^c)} = \dfrac {|\partial A|} {{\vol(A)}}=h_2(G)
\end{split}
\end{equation*}

Now, let us consider a set $B\subseteq V$ such that  $B\not\in\{\emptyset, A, A^c, V\}$ and 
\begin{equation}
\label{cond_3_part}
\vol(B)\le\vol(B^c),\quad \frac{|\partial B|}{\vol(B)},\frac{|\partial B^c|}{\vol(B^c)}\le \frac{|\partial (A\cap B^c)|+|\partial (A^c\cap B)|}{\vol(A\cap B^c)+\vol(A^c \cap B)}.
\end{equation}
The triple $\{B, A\cap B^c, A^c\cap B^c\}$ is a partition of $V$ and therefore, by definition, we have
\[
h_3(G)\leq \max\left\{ \dfrac {|\partial B|} {{\vol(B)}},  \dfrac {|\partial (A\cap B^c)|} {{\vol(A \cap B^c)}}, \dfrac {|\partial (A^c\cap B^c)|} {{\vol(A^c\cap B^c)}}\right\} = \dfrac {|\partial B|} {{\vol(B)}},
\]
where the last equality holds up to rename the sets. Hence, we have
\begin{equation}
\label{optimization_c1c2}
\begin{split}
&\max_{c_1,c_2\in\R} \dfrac{\sum_{\substack{i,j\in V\\ j\sim i}}|(  1_B)_i-(  1_B)_j-c_2((1_A)_i-(  1_A)_j)|}{\sum_{i\in V}d_i|(  1_B)_i-c_1(  1_V)_i-c_2(  1_A)_i|}\\
&=\max_{\substack{c_1,c_2\ge 0\\c_1+c_2\leq 1}}\dfrac{\substack{E(A\cap B, A\cap B^c)+c_2E(A\cap B, A^c\cap B)+(1-c_2)E(A\cap B, A^c\cap B^c)\\+(1+c_2)E(A\cap B^c, A^c\cap B)+c_2E(A\cap B^c, A^c\cap B^c)+E(A^c\cap B^c, A^c\cap B)}}{\substack{(1-c_1-c_2)\vol (A\cap B)+(c_1+c_2)\vol (A\cap B^c)\\+(1-c_1)\vol (A^c\cap B)+c_1\vol (A^c\cap B^c)}}
= \dfrac {|\partial B|} {{\vol(B)}}\ge h_3(G)
\end{split}
\end{equation}
We say that ${\bf y}_3={\bf 1}_{\bar B}$, where $\bar B$ is the set achieving the minimum in the first term of \eqref{optimization_c1c2} among sets verifying \eqref{cond_3_part}, that is
\begin{equation*}
\begin{split}
&\max_{c_1, c_2 \in\R}  \dfrac{\sum_{\substack{i,j\in V\\ j\sim i}}|(y_3)_i-(y_3)_j-c_2((  y_2)_i-(y_2)_j)|}{ \sum_{i\in V}d_i|g_i-c_1-c_2(  1_A)_i|}\\
&=\min_{\substack{B\not\in\{\emptyset, A, A^c, V\}\\
\vol(B)\le\vol(B^c)\\ \frac{|\partial B|}{\vol(B)},\frac{|\partial B^c|}{\vol(B^c)}\le \frac{|\partial (A\cap B^c)|+|\partial (A^c\cap B)|}{\vol(A\cap B^c)+\vol(A^c \cap B)}}}\max_{c_1,c_2\in\R} \dfrac{\sum_{\substack{i,j\in V\\ j\sim i}}|(  1_B)_i-(  1_B)_j-c_2((1_A)_i-(  1_A)_j)|}{\sum_{i\in V}d_i|(  1_B)_i-c_1(  1_V)_i-c_2(  1_A)_i|}\ge h_3(G).
\end{split}
\end{equation*}
\end{proof}
\begin{rem}
We stress that the proof of the previous result implies that the inequality stated for the third Cheeger constant holds as an equality if two of the three sets realizing $h_3(G)$ are entirely contained in $A$ or in $A^c$.
\end{rem}

\subsection{The pseudo-orthogonality}
In this Section, we introduce the concept of pseudo-orthogonality and we use it to study the critical points of the functional $\hat I$. Throughout this Section, for any couple of matrices $A=(a_{ij}),B=(b_{ij})\in\R^{n\times n}$, we denote the matrix product
\[
C=AB\quad\text{where}\quad c_{ij}=\sum_{h=1}^n a_{ih}b_{hj}\quad\forall i,j=1,...,n,
\]
and the Hadamard product 
\[
C=A\odot B\quad\text{where}\quad  c_{ij}=a_{ij}b_{ij}\quad\forall i,j=1,...,n.
\]
Moreover, we denote $W=(w_{ij})\in\R^{n\times n}$ the weight matrix, that is the symmetric matrix defined such that $w_{ij}$ is equal to 1 or 0 if $i \sim j$ or $i \not\sim j$, respectively. \begin{prop} \label{ott_funz}
Let $G=(V,E)$ be a graph and $\mathbf g$ a vector of $\R^n$. Let ${\bf 1}_V$ and ${\bf f}_2$ the first and the second eigenfunctions of the graph $1$-Laplacian eigenvalue problem \eqref{eigen}. Then the following holds.
\begin{enumerate}
\item The critical points of the function $c\in\R\mapsto ||{\bf g}-c{\bf 1}_V||_w$ are achieved for $\bar c$ such that $0\in\langle D\Sgn ({\bf g}-\bar c{\bf  1}_V), {\bf 1}_V\rangle$. 
Moreover, for any $\bar{\bf g}$ with $0\in \langle D \Sgn\bar{\bf g},{\bf 1}_V\rangle$, we have $|| \bar{\bf g} ||_w=\min_{c\in\R}||\bar{\bf g}-c {\bf 1}_V||_w$.
\item The critical points of the function $(c_1, c_2)\in\R^2\mapsto ||{\bf g}-c_1{\bf 1}_V-c_2{\bf f}_2||_w$ are achieved for $(\bar c_1,\bar c_2)$ such that $0\in\langle D\Sgn ({\bf g}-\bar c_1{\bf  1}_V-\bar c_2 {\bf f}_2), {\bf 1}_V\rangle$ and $0\in\langle D\Sgn ({\bf g}-\bar c_1{\bf  1}_V-\bar c_2 {\bf f}_2), {\bf f}_2\rangle$. Moreover, for any $\bar{\bf g}$ with $0\in\langle D\Sgn (\bar{\bf g}), {\bf 1}_V\rangle$ and $0\in\langle D\Sgn (\bar{\bf g}), {\bf f}_2\rangle$, we have  $||\bar {\bf g} ||_w=\min_{c_1, c_2\in\R}||{\bf g}-c_1{\bf 1}_V-c_2{\bf f}_2||_w$. 
\item The critical points of the function $c_2\in\R\mapsto I({\bf g}-c_2{\bf f}_2)$ are achieved for $\bar c_2$ such that $0\in\langle (w_{ij})\odot (\Sgn ((g_i-g_j)- \bar c_2 (({ f_2})_i-({ f_2})_j))\odot (( f_2)_i-(f_2)_j){\bf 1}_V,{\bf 1}_V\rangle$. 
Moreover, for any $\bar{\bf g}$ with  $0\in\langle (w_{ij})\odot (\Sgn (\bar g_i-\bar g_j))\odot (( f_2)_i-(f_2)_j){\bf 1}_V,{\bf 1}_V\rangle$, we have $I(\bar{\bf g})=\max_{ c_2\in\R} I(\bar{\bf g}-c_2{\bf f}_2)$.
\item The critical points of the function $(c_1,c_2)\in\R\mapsto \hat I({\bf g}-c_1 {\bf 1}_V-c_2{\bf f}_2)$ are achieved for $(\bar c_1,\bar c_2)$ such that $0\in\langle D\Sgn ({\bf g}-\bar c_1{\bf  1}_V-\bar c_2 {\bf f}_2), {\bf 1}_V\rangle$  and $0\in\langle (w_{ij})\odot (\Sgn ((g_i-g_j)- \bar c_2 (({ f_2})_i-({ f_2})_j))\odot (( f_2)_i-(f_2)_j){\bf 1}_V,{\bf 1}_V\rangle||{\bf g}-\bar c_1{\bf 1}_V-\bar c_2 {\bf f}_2||_w-I({\bf g}-\bar c_1{\bf  1}_V-\bar c_2 {\bf f}_2) \langle D\Sgn ({\bf g}-\bar c_1{\bf 1}_V-\bar c_2 {\bf f}_2), {\bf f}_2\rangle$. 
Moreover, for any $\bar{\bf g}$ with $0\in \langle D \Sgn (\bar{\bf g}),{\bf 1}_V\rangle$ and $0\in\langle (w_{ij})\odot (\Sgn (\bar g_i-\bar g_j))\odot (( f_2)_i-(f_2)_j){\bf 1}_V,{\bf 1}_V\rangle||{\bf g}||_w-I(\bar {\bf g})\langle D\Sgn (\bar{\bf g}), {\bf f}_2\rangle$, we have $\hat I(\bar{\bf g})=\max_{c_1, c_2\in\R} \hat I(\bar{\bf g}-c_1{\bf 1}_V-c_2{\bf f}_2)$.
\end{enumerate}
\end{prop}
\begin{proof}{\it (1)}
The critical points $\bar c$ are such that:
\[
0\in \sum_{i\in V}d_i\Sgn\left(g_i-\bar c\right)
=\langle D\Sgn ({\bf g}-\bar c{\bf  1}_V), {\bf 1}_V\rangle.
\]
In particular, we are saying that $\bar c$ is the weighted median of ${\bf g}$. 

{\it (2)} The critical points $(\bar c_1, \bar c_2)$ are such that:
\begin{align*}
0\in{\sum_{i\in V}d_i \Sgn \left(g_i-\bar c_1-{\bar c_2}({ f_2})_i\right)}
=\langle D\Sgn ({\bf g}-\bar c_1 {\bf 1}_V-\bar c_2 {\bf f}_2), {\bf 1}_V\rangle;
\end{align*}
\begin{align*}
0&\in{\sum_{i\in V}d_i\Sgn \left(g_i- {\bar c_1}-{\bar c_2}({ f_2})_i\right)}({ f_2})_i
=\langle D\Sgn ({\bf g}-\bar c_1{\bf  1}_V-\bar c_2 {\bf f}_2), {\bf f}_2\rangle.
\end{align*}

{\it (3)} The critical points $ \bar c_2$ are such that:
\begin{equation*}
\begin{split}
0\in\sum_{\substack{i,j\in V\\ i\sim j}}&\Sgn (g_i-g_j-\bar c_2 (({ f_2})_i-({ f_2})_j))\cdot (( f_2)_i-( f_2)_j) \\
=&\langle (w_{ij})\odot (\Sgn ((g_i-g_j)-\bar c_2 (({ f_2})_i-({ f_2})_j))\odot (( f_2)_i-(f_2)_j){\bf 1}_V,{\bf 1}_V\rangle
\end{split}
\end{equation*}

{\it (4)} The critical points $(\bar c_1, \bar c_2)$ are such that:
 \begin{equation*}
0\in-\dfrac{I({\bf f}-\bar c_1{\bf  1}_V-\bar c_2 {\bf f}_2)}{||{\bf f}-\bar c_1{\bf 1}_V-\bar c_2 {\bf f}_2||_w^2}\langle D\Sgn ({\bf f}-\bar c_1{\bf  1}_V-\bar c_2 {\bf f}_2), {\bf 1}_V\rangle;
\end{equation*}
 \begin{equation*}
    \begin{split}
0\in&\frac{\sum_{\substack{i,j\in V\\ i\sim j}}\Sgn (f_i-f_j-\bar c_2 (({\hat f_2})_i-({ f_2})_j))\cdot (( f_2)_i-(f_2)_j) \cdot ||{\bf f}-\bar c_1{\bf 1}_V-\bar c_2 {\bf f}_2||_w}{||{\bf f}-\bar c_1{\bf 1}_V-\bar c_2 {\bf f}_2||_w^2}\\
&-\frac{I({\bf f}-\bar c_1{\bf  1}_V-\bar c_2 {\bf f}_2) \langle D\Sgn ({\bf f}-\bar c_1{\bf 1}_V-\bar c_2 {\bf f}_2), {\bf f}_2\rangle}{||{\bf f}-\bar c_1{\bf 1}_V-\bar c_2 {\bf f}_2||_w^2}\\
=&\frac{\langle (w_{ij})\odot (\Sgn ((g_i-g_j)-\bar c_2 (({ f_2})_i-({ f_2})_j))\odot (( f_2)_i-(f_2)_j){\bf 1}_V,{\bf 1}_V\rangle ||{\bf f}-\bar c_1{\bf 1}_V-\bar c_2 {\bf f}_2||_w}{||{\bf f}-\bar c_1{\bf 1}_V-\bar c_2 {\bf f}_2||_w^2}\\
&-\frac{I({\bf f}-\bar c_1{\bf  1}_V-\bar c_2 {\bf f}_2) \langle D\Sgn ({\bf f}-\bar c_1{\bf 1}_V-\bar c_2 {\bf f}_2), {\bf f}_2\rangle}{||{\bf f}-\bar c_1{\bf 1}_V-\bar c_2 {\bf f}_2||_w^2}.
    \end{split}
\end{equation*}
\end{proof}
Proposition \ref{ott_funz} leads to the following inductive definition of pseudo-orthogonality.
\begin{defn}\label{ort_def}
Let $G=(V,E)$ be a graph and $\hat{\bf g}_k$ the vector realizing $m_k(G)$ as defined in \eqref{ort_value}. Since we know that ${\bf g}_1={\bf 1}_V$, we say a vector ${\bf g}$ is pseudo-orthogonal to ${\bf g}_1$ and we denote
\[
{\bf g}\perp_p{\bf g}_1\quad\Longleftrightarrow\quad0\in \langle D \Sgn ({\bf g}), {\bf 1}_V\rangle,
\]
that is when ${\bf g}$ has zero weighted median.

Since we know that ${\bf g}_2={\bf f}_2$, we say a vector ${\bf g}$ is pseudo-orthogonal to ${\bf g}_2$ and we denote
\begin{equation*}
{\bf g}\perp_p{\bf g}_2\Longleftrightarrow   
\begin{cases}
0\in \langle D \Sgn ({\bf g}), {\bf 1}_V\rangle,\\
0\in\langle (w_{ij})\odot (\Sgn (g_i- g_j))\odot (( f_2)_i-(f_2)_j){\bf 1}_V,{\bf 1}_V\rangle-\hat I({\bf g})\langle D\Sgn ({\bf g}), {\bf f}_2\rangle.
\end{cases}
\end{equation*}
\end{defn}
\begin{rem}
Let us observe that the following conditions
\begin{equation}\label{g_2_ort}
    \begin{split}
&0\in \langle D \Sgn ({\bf g}), {\bf 1}_V\rangle,\\
&0\in \langle D \Sgn ({\bf g}), {\bf f}_2\rangle,\\
&0\in\langle (w_{ij})\odot (\Sgn (g_i- g_j))\odot (( f_2)_i-(f_2)_j){\bf 1}_V,{\bf 1}_V\rangle
    \end{split}
\end{equation}
imply that ${\bf g}$ is pseudo-orthogonal to ${\bf g}_1$ and ${\bf g}_2$. This observation could be further investigate to state that, for any $k\in\N$, the following conditions
\begin{equation*}
    \begin{split}
&0\in \langle D \Sgn ({\bf g}), {\bf g}_m\rangle\quad\forall m=1,...,k,\\
&0\in\langle (w_{i,j})\odot (
\Sgn ( g_i-g_j))\odot (( g_m)_i-( g_m)_j){\bf 1}_V,{\bf 1}_V\rangle \quad\forall m=2,...,k.
    \end{split}
\end{equation*}
imply that ${\bf g}$ is pseudo-orthogonal to ${\bf g}_1$, ..., ${\bf g}_k$.
\end{rem}

\begin{prop}
Let $G=(V,E)$ be a graph. Then $m_k(G)$ is an eigenvalue of \eqref{graph_1} for any $k\geq 2$.
\end{prop}
\begin{proof}
The result for $k=2$ easily follows by Proposition \ref{2_eq_chain}. For $k\geq 3$, let us consider $\hat{\bf g}_k\in X$ (one of) the vector realizing $m_k(G)$. Then, since $m_k(G)$ is a critical point for the functional $I$ on the subset of vectors pseudo-orthogonal to the vectors realizing the previous constants $m_1(G), ..., m_{k-1}(G)$, then
\[
0\in\frac{d}{dt}\left(\sum_{i\sim j}|(g_k)_i-(g_k)_j|+t(u_i-u_j)-m_k(G)||{\bf g}_k+t{\bf u}||_w\right)|_{t=0}
\]
for any ${\bf u}\in\R^n$. By choosing any ${\bf u}=(u_1,0,...,0)$, we have
\[
0\in \sum_{1\sim j}\Sgn((g_k)_1-(g_k)_j)u_1-m_k(G)d_i \Sgn((g_k)_1)u_1.
\]
This means that
\[
0\in \sum_{1\sim j}\Sgn((g_k)_1-(g_k)_j)-m_k(G)d_i \Sgn((g_k)_1),
\]
that gives the first component of eigenpair as in \eqref{eigen}. Analogously the other $n-1$ components are obtained and the result follows.
\end{proof}

\subsection{Proof of the main results} Now we consider a characterization through the span of the first and second eigenfunctions of the $1$-Laplacian.
\begin{prop}\label{span_h2h3}
Let $G=(V,E)$ be a graph. Then: 
\begin{equation*}
\begin{split}
h_2(G)&=\min_{{\bf g}\not\in< {\bf 1}_V>}\max_{c\in\R} \dfrac{\sum_{\substack{i,j\in V\\ j\sim i}}|g_i-g_j|}{\sum_{i\in V}d_i|g_i-c|},\\
h_3(G)&\leq\min_{{\bf g}\not\in< {\bf 1}_V,{\bf f}_2>}\max_{c_1, c_2 \in\R} \dfrac{\sum_{\substack{i,j\in V\\ j\sim i}}{|g_i-g_j-c_2((f_2)_i-(f_2)_j)|}}{\sum_{i\in V}d_i|g_i-c_1-c_2(f_2)_i|}.
\end{split}
\end{equation*}
\end{prop}
\begin{proof}
Since ${\bf y}_2$ 
in Proposition \ref{espl_h2h3} is not in $<{\bf 1}_V
>$, we only need to prove that
\begin{equation*}
\begin{split}
h_2(G)&\leq\min_{{\bf g}\not\in< {\bf 1}_V>}\max_{c\in\R} \dfrac{\sum_{\substack{i,j\in V\\ j\sim i}}|g_i-g_j|}{\sum_{i\in V}d_i|g_i-c|},\\
h_3(G)&\leq\min_{{\bf g}\not\in< {\bf 1}_V,{\bf f}_2>}\max_{c_1, c_2 \in\R} \dfrac{\sum_{\substack{i,j\in V\\ j\sim i}}{|g_i-g_j-c_2((f_2)_i-(f_2)_j)|}}{\sum_{i\in V}d_i|g_i-c_1-c_2(f_2)_i|}.
\end{split}
\end{equation*}
For any ${\bf g}\not\in<{\bf 1}_V>$, let us fix $\bar c$ such that $0\in\langle D\Sgn ({\bf g}-\bar c{\bf 1}_V), {\bf 1}_V\rangle$
.
Then, for any $\sigma\in\R$, we consider a function counting the edges between the superlevel set and the sublevel set of ${\bf g}-\bar c{\bf 1}_V$:
\[
G(\sigma)=|\{(i,j)\in E\ | \ g_i-\bar c \leq \sigma<g_j-\bar c\}|.
\]
Therefore, we have
\begin{equation*}
\begin{split}
&\dfrac{\sum_{\substack{i,j\in V\\ j\sim i}}|g_i-g_j|}{\sum_{i\in V}d_i|g_i-\bar c|} =\dfrac{\dint_{-\infty}^{+\infty}G(\sigma)\ \text{d}\sigma}{\sum_{i\in V}d_i|g_i-\bar c|}\\
&=\dfrac{\dint_{-\infty}^{0}\frac{G(\sigma)}{\sum_{g_i-\bar c<\sigma}d_i}\sum_{g_i-\bar c<\sigma}d_i\ \text{d}\sigma+\dint_{0}^{+\infty}\frac{G(\sigma)}{\sum_{g_i-\bar c>\sigma}d_i}\sum_{g_i-\bar c>\sigma}d_i\ \text{d}\sigma}{\sum_{i\in V}d_i|g_i-\bar c|}\\
&\geq h_2(G)\dfrac{\dint_{-\infty}^{0}\sum_{g_i-\bar c<\sigma}d_i\ \text{d}\sigma+\dint_{0}^{+\infty}\sum_{g_i-\bar c>\sigma}d_i\ \text{d}\sigma}{\sum_{i\in V}d_i|g_i-\bar c|}=h_2(G).
\end{split}
\end{equation*}
Hence the conclusion for the second Cheeger constant follows by passing to the supremum for any real constant and to the infimum for any ${\bf g}\not\in< {\bf 1}_V>$.

Now, for any ${\bf g}\not\in< {\bf 1}_V,{\bf f}_2>$, let us fix $\bar c_1,\bar c_2$ such that $0\in\langle D\Sgn ({\bf g}-\bar c_1{\bf  1}_V-\bar c_2{\bf f}_2, {\bf 1}_V\rangle$, $0\in\langle (w_{ij})\odot (\Sgn (\bar g_i-\bar g_j))\odot (( f_2)_i-(f_2)_j){\bf 1}_V,{\bf 1}_V\rangle||{\bf g}||_w-I({\bf g})\langle D\Sgn ({\bf g}), {\bf f}_2\rangle$.
Then, for any $\sigma\in\R$, we consider a function counting the edges between the superlevel set and the sublevel set of ${\bf g}-\bar c_1 {\bf 1}_V -\bar c_2{\bf  f}_2$:
\begin{equation*}
G(\sigma)=\left|\left\{(i,j)\in E\ :
\begin{aligned} \ g_i-\bar c_1-\bar c_2 (f_2)_i\leq \sigma<g_j-\bar c_1-\bar c_2 (f_2)_i
\end{aligned}
\right\}\right|.
\end{equation*}
Therefore, we have
\begin{equation*}
\begin{split}
&\dfrac{\sum_{\substack{i,j\in V\\ j\sim i}}{|g_i-g_j-\bar c_2((f_2)_i-(f_2)_j)|}}{\sum_{i\in V}d_i|g_i-\bar c_1-\bar c_2(f_2)_i|} =\dfrac{\dint_{-\infty}^{+\infty}G(\sigma)\ \text{d}\sigma}{\sum_{i\in V}d_i|g_i-\bar c_1-\bar c_2 (f_2)_i|}\\
&=\dfrac{\dint_{-\infty}^{0}\frac{G(\sigma)}{\sum_{g_i-\bar c_1-\bar c_2 (f_2)_i<\sigma}d_i}\sum_{g_i-\bar c_1-\bar c_2 (f_2)_i<\sigma}d_i\ \text{d}\sigma+\dint_{0}^{+\infty}\frac{G(\sigma)}{\sum_{g_i-\bar c_1-\bar c_2 (f_2)_i>\sigma}d_i}\sum_{g_i-\bar c_1-\bar c_2 (f_2)_i>\sigma}d_i\ \text{d}\sigma}{\sum_{i\in V}d_i|g_i-\bar c_1-\bar c_2 (f_2)_i|}\\
&\geq h_3(G)\dfrac{\dint_{-\infty}^{0}\sum_{g_i-\bar c_1-\bar c_2 (f_2)_i<\sigma}d_i\ \text{d}\sigma+\dint_{0}^{+\infty}\sum_{g_i-\bar c_1-\bar c_2 (f_2)_i>\sigma}d_i\ \text{d}\sigma}{\sum_{i\in V}d_i|g_i-\bar c_1-\bar c_2 (f_2)_i|}=h_3(G).
\end{split}
\end{equation*}
Hence the conclusion for the third Cheeger constant follows by passing to the supremum for any couple of real constants and to the infimum for any ${\bf g}\not\in< {\bf 1}_V,{\bf f}_2>$.
\end{proof}

Therefore, by Propositions \ref{ott_funz}, \ref{espl_h2h3} and \ref{span_h2h3}, we prove the following.
\begin{thm}\label{hNCC_thm}
Let $G=(V,E)$ be a graph, then:
\begin{equation*}
    \begin{split}
    (i) &\quad    h_2(G)=m_2(G),\\
    (ii) &\quad     m_3(G)=\min_{{\bf g}\not\in< {\bf 1}_V,{\bf f}_2>}\max_{c_1, c_2 \in\R} \dfrac{\sum_{\substack{i,j\in V\\ j\sim i}}{|g_i-g_j-c_2((f_2)_i-(f_2)_j)|}}{\sum_{i\in V}d_i|g_i-c_1-c_2(f_2)_i|},\\
     (iii) &\quad     h_3(G)\leq m_3(G),
    \end{split}
\end{equation*}

\end{thm}
\begin{proof}
To prove {\it (i)}, we need to prove two inequalities.
\begin{itemize}
    \item 
Firstly, we prove $h_2(G)\ge m_2(G)$. From Proposition \ref{espl_h2h3}, we know there exists ${\bf y}_2$ such that
\[
h_2(G)=\max_{c\in\R} \dfrac{\sum_{\substack{i,j\in V\\ j\sim i}}|(y_2)_i-(y_2)_j|}{\sum_{i\in V}d_i|(y_2)_i-c|}.
\]
Let us set $\bar c$ such that $0\in\langle D\Sgn ({\bf y}_2-\bar c{\bf  1}_V), {\bf 1}_V\rangle$ and hence ${\bf z}_2 ={\bf y}_2-\bar c {\bf 1}_V$. Therefore ${\bf z}_2\perp_p {\bf 1}_V$, and we have:
\begin{equation*}
    \begin{split}
h_2(G)&=\max_{c\in\R} \dfrac{\sum_{\substack{i,j\in V\\ j\sim i}}|(y_2)_i-(y_2)_j|}{\sum_{i\in V}d_i|(y_2)_i-c|}\ge \dfrac{\sum_{\substack{i,j\in V\\ j\sim i}}|(y_2)_i-(y_2)_j|}{\sum_{i\in V}d_i|(y_2)_i-\bar c|}\\
&\ge \dfrac{\sum_{\substack{i,j\in V\\ j\sim i}}|(z_2)_i-(z_2)_j|}{\sum_{i\in V}d_i|(z_2)_i|}\geq \min_{{\bf z}\perp_p {\bf 1}_V} \dfrac{\sum_{\substack{i,j\in V\\ j\sim i}}|z_i-z_j|}{\sum_{i\in V}d_i|z_i|} =m_2(G),
    \end{split}
\end{equation*}
where $m_2(G)$ is defined in \eqref{ort_value}.
\item Now, we prove that $m_2(G)\ge h_2(G)$. Let us denote by $\hat{\bf g}_2$ a vector in $X$ pseudo-orthogonal to ${\bf 1}_V$ such that
\[
m_2(G)=\sum_{\substack{i,j\in V\\ j\sim i}}|(\hat g_2)_i-(\hat g_2)_j|.
\]
Let us set $\bar c$ such that $0\in\langle D\Sgn (\hat {\bf g}_2-\bar c{\bf  1}_V), {\bf 1}_V\rangle$, then by Propositions \ref{ott_funz} (1) and \ref{span_h2h3}, we have
\begin{equation*}
\begin{split}
m_2(G)&=\sum_{\substack{i,j\in V\\ j\sim i}}|(\hat g_2)_i-(\hat g_2)_j|=\max_{c\in\R} \dfrac{\sum_{\substack{i,j\in V\\ j\sim i}}|(\hat g_2)_i- (\hat  g_2)_j|}{\sum_{i\in V}d_i|(\hat  g_2)_i-c|}\\
&\ge\inf_{{\bf g}\not\in\langle \hat{\bf 1}_V\rangle}\sup_{c\in\R} \dfrac{\sum_{\substack{i,j\in V\\ j\sim i}}|g_i-g_j|}{\sum_{i\in V}d_i|g_i-c|}=h_2(G).
\end{split}
\end{equation*}
\end{itemize}
To prove {\it (ii)}, we need to show two inequalities.
 \begin{itemize}
\item  From Proposition \ref{espl_h2h3}, we know that there exist ${\bf y}_2$ and ${\bf y}_3$ such that
\[
\max_{c_1, c_2 \in\R} \dfrac{\sum_{\substack{i,j\in V\\ j\sim i}}{|(y_3)_i-(y_3)_j-c_2((y_2)_i-(y_2)_j)|}}{\sum_{i\in V}d_i|(y_3)_i-c_1-c_2(y_2)_i|}.
\]
Let us set $\bar c_1$ and $\bar c_2$ such that $0\in\langle D\Sgn ({\bf y}_3-\bar c_1{\bf  1}_V-\bar c_2 {\bf y}_2), {\bf 1}_V\rangle$ and $0\in\langle (w_{ij})\odot (\Sgn (((y_3)_i-(y_3)_j)- \bar c_2 (({ y_2})_i-({y_2})_j))\odot (( y_2)_i-(y_2)_j){\bf 1}_V,{\bf 1}_V\rangle||{\bf y}_3-\bar c_1{\bf 1}_V-\bar c_2 {\bf y}_2||_w-I({\bf y}_3-\bar c_1{\bf  1}_V-\bar c_2 {\bf y}_2) \langle D\Sgn ({\bf y}_3-\bar c_1{\bf 1}_V-\bar c_2 {\bf y}_2), {\bf y}_2\rangle$ and hence
${\bf z}_3 ={\bf y}_3-\bar c_1 {\bf 1}_V-\bar c_2 {\bf y}_2$. Therefore ${\bf z}_3\perp_p {\bf y}_2$ and ${\bf z}_3\perp_p {\bf 1}_V$, and we have:
\begin{equation*}
\begin{split}
    &\max_{c_1, c_2 \in\R} \dfrac{\sum_{\substack{i,j\in V\\ j\sim i}}{|(y_3)_i-(y_3)_j-c_2((y_2)_i-(y_2)_j)|}}{\sum_{i\in V}d_i|(y_3)_i-c_1-c_2(y_2)_i|}\\
    &\ge\dfrac{\sum_{\substack{i,j\in V\\ j\sim i}}{|(y_3)_i-(y_3)_j-\bar c_2((y_2)_i-(y_2)_j)|}}{\sum_{i\in V}d_i|(y_3)_i- \bar c_1-\bar c_2(y_2)_i|}\\
    &\ge         \dfrac{\sum_{\substack{i,j\in V\\ j\sim i}}{|(z_3)_i-(z_3)_j|}}{\sum_{i\in V}d_i|(z_3)_i|}\ge\min_{\substack{{\bf z} \perp_p \hat{\bf 1}_V\\ {\bf z}\perp_p \hat{\bf f}_2}} \dfrac{\sum_{\substack{i,j\in V\\ j\sim i}}{|z_i-z_j|}}{\sum_{i\in V}d_i|z_i|}=m_3(G).
\end{split}
\end{equation*}
\item Now, let us denote by $\hat {\bf g}_3$ a vector in $X$ pseudo-orthogonal to ${\bf 1}_V$ and ${\bf f}_2$ such that
\[
m_3(G)=\sum_{\substack{i,j\in V\\ j\sim i}}|(\hat g_3)i-(\hat g_3)_j|.
\]
Let us set $\bar c_1$ and $\bar c_2$ such that $0\in\langle D\Sgn (\hat {\bf g}_3-\bar c_1{\bf  1}_V-\bar c_2 {\bf f}_2), {\bf 1}_V\rangle$ and $0\in\langle (w_{ij})\odot (\Sgn (((\hat g_3)_i-(\hat g_3)_j)- \bar c_2 (({ f_2})_i-({ f_2})_j))\odot (( f_2)_i-(f_2)_j){\bf 1}_V,{\bf 1}_V\rangle||\hat{\bf g}_3-\bar c_1{\bf 1}_V-\bar c_2 {\bf f}_2||_w-I(\hat{\bf g}_3-\bar c_1{\bf 1}_V-\bar c_2 {\bf f}_2) \langle D\Sgn (\hat{\bf g}_3-\bar c_1{\bf 1}_V-\bar c_2 {\bf f}_2), {\bf f}_2\rangle$. Then by Propositions \ref{ott_funz}(4) and \ref{span_h2h3}, we have
\begin{equation*}
\begin{split}
m_3(G)
&=\sum_{\substack{i,j\in V\\ j\sim i}}|(\hat g_3)i-(\hat g_3)_j|\\
& =\sup_{c_1, c_2 \in\R} \dfrac{\sum_{\substack{i,j\in V\\ j\sim i}}{|(\hat g_3)_i-(\hat g_3)_j- c_2(( f_2)_i-( f_2)_j)|}}{\sum_{i\in V}d_i|(\hat g_3)_i- c_1- c_2( f_2)_i|}\\
&\ge\min_{{\bf g}\not\in\langle {\bf 1}_V,{\bf f}_2\rangle}\max_{c_1, c_2 \in\R} \dfrac{\sum_{\substack{i,j\in V\\ j\sim i}}{|g_i-g_j-c_2((f_2)_i-(f_2)_j)|}}{\sum_{i\in V}d_i|g_i-c_1-c_2(f_2)_i|}.
\end{split}
\end{equation*}
 \end{itemize}
 Finally, the claim {\it (iii)} easily follows from Proposition \ref{span_h2h3} and the previous point:
 \[
 h_3(G)\leq\min_{{\bf g}\not\in< {\bf 1}_V,{\bf f}_2>}\max_{c_1, c_2 \in\R} \dfrac{\sum_{\substack{i,j\in V\\ j\sim i}}{|g_i-g_j-c_2((f_2)_i-(f_2)_j)|}}{\sum_{i\in V}d_i|g_i-c_1-c_2(f_2)_i|}=m_3(G).
 \]
\end{proof}

\begin{proof}{\it Theorem \ref{realizing}.}
The desired chain of inequalities \eqref{3chain} follows by Theorems \ref{ck_rhok} and \ref{hNCC_thm}.
 \end{proof}

\begin{rem}
\label{k_ort_rem}
Following the ideas exposed in this paper, it would be hopeful to characterize the third Cheeger constant as the minimum of functional \eqref{I} among vectors pseudo-orthogonal to ${\bf g}_1={\bf 1}_V$ and ${\bf g}_2={\bf f}_2$. Furthermore, it would be reasonable to generalize these results to $k$-Cheeger constant, for $k>3$. More precisely, we expect that the $k$-th Cheeger constant is the minimum of \eqref{I} among vectors $\hat{\bf g}$ such that 
\begin{equation}
    \label{g_k_ort}
\hat{\bf g}\perp_p {\bf g}_1,\quad \hat{\bf g} \perp_p {\bf g}_2,\quad...,\quad\hat{\bf g}\perp_p {\bf g}_{k-1}.
\end{equation}
\end{rem}

\section{Application of the Inverse Power Method to Spectral Data Clustering}\label{AotIPMtSPD}
We perform the $1$-Spectral Clustering based on the inverse power method.
The inverse power method (IPM) is a standard technique to obtain the smallest eigenvalue of a positive semi-definite symmetric matrix $A$ based on the following iterative scheme:
\[
A f^{k+1}=f^k\quad k\in\N,
\] 
{trasformed in the optimization problem}:
\[
f^{k+1}=\arg \min_u \frac 12 (u, Au)-(u,f^k)\quad k\in\N.
\]
The IPM can be extended to nonlinear cases as in \cite{BHb}.

Before explain how our algorithm works, we give the definition of Cheeger constants in a more treatable way for numeric applications. 
The main input data we need to perform the algorithm we are presenting, is the weight matrix $W=(w_{ij})_{i,j=1}^n$, that is a symmetric $n\times n$ matrix defined such that $w_{ij}$ is equal to 1 or 0 if $i \sim j$ or $i \not\sim j$, respectively. For any subset $A\subseteq V$ we call the cut of $A$ the quanity
\[
\cut(A,A^c)=\sum_{i\in A, j\in A^c} w_{ij}.
\]
Moreover, given $A,B,C\subseteq V$, we denote the normalized $2$-Cheeger cut and the normalized $3$-Cheeger cut the quantities
\begin{equation*}
\begin{split}
&\NCC_2(A,B)=\max\left\{\dfrac{\cut(A,A^c)}{\vol(A)},\dfrac{\cut(B,B^c)}{\vol(B)}\right\},\\
&\NCC_3(A,B,C)=\max\left\{\dfrac{\cut(A,A^c)}{\vol(A)},\dfrac{\cut(B,B^c)}{\vol(B)},\dfrac{\cut(C,C^c)}{\vol(C)}\right\}.
\end{split}
\end{equation*}  
Furthermore, we call the second and the third optimal normalized Cheeger cut (that are the second and the third Cheeger constant, respectively) the quantities:
\begin{equation*}
\begin{split}
h_2(G):=\inf_{A,B\subseteq V} NCC_2 (A, B),\quad
h_3(G):=\inf_{A,B,C\subseteq V} NCC_3 (A,B, C).
\end{split}
\end{equation*}

The algorithm we propose is based on a transformation of the graph Cheeger problem \eqref{cheeger_constant} into a problem of optimizing the functional \eqref{I}. 
Therefore the vector realizing the third Cheeger constant is characterized as in \eqref{g_2_ort}. 

We modify the algorithm that has been proposed in \cite{BHb}. Particularly, ${\bf f}_2$ is computed by an iteration in which each times the weighted median is subtracted. Similarly, the vector realizing the third Cheeger constant is obtained by an iterative process in which each time the weighted median is subtracted and the obtained vector is worked by the routine $PseudoOrt$. Specifically, this routine realizes the second pseudo orthogonality condition in Definition \ref{ort_def}: from a starting vector is subtracted $\lambda \hat{\bf f}_2$, for suitable real constant $\lambda$.
\begin{center}
{\bf Algorithm to compute the second eigenvector}\\
\vspace{-0.4cm}
\begin{tabular}{l}
  \\ \hline
 {\bf Initialization} with $\mathbf f_0$ non costant vector such that $\median({\bf f_0})=0$ and $||{\bf f_0}||_1=1$  \\ \hline
 {\bf Repeat}\qquad\qquad\qquad \\ 
1. $g^{k+1}=\argmin_{||f||_2^2\leq 1}\left\{\frac 12 \sum_{i, j =1}^n\omega_{ij}|f_i-f_j|-\mu^k(f\cdot v^k)\right\}$\\ 
2. $ f^{k+1}=g^{k+1}-\median(g^{k+1})$\\ 
3. $v_i^{k+1}=\begin{cases}\sign(f_i^{k+1}), \qquad\qquad\text{if}\ f_i^{k+1}\neq 0,\\ -\dfrac{|f_+^{k+1}|-|f_-^{k+1}|}{|f_0^{k+1}|},\quad\text{if} \ f_i^{k+1}=0. \end{cases}$\\
4. $\mu^{k+1}=I(f^{k+1})$\\
{\bf Until} $\dfrac{|\mu^{k+1}-\mu^k|}{\mu^k}<\eps$\\ \hline
\end{tabular}
\end{center}

\begin{center}
{\bf Algorithm to compute the third eigenvector}\\
\vspace{-0.4cm}
\begin{tabular}{l}
  \\ \hline
 {\bf Initialization} with ${\bf f_0}$ non costant vector such that $\median({\bf f_0})=0$ and $(\sign({\bf f_0})\cdot {\bf u_2})=0$  \\ \hline
 {\bf Repeat}\qquad\qquad\qquad \\ 
1. $g^{k+1}=\argmin_{||f||_2^2\leq 1}\left\{\frac 12 \sum_{i, j =1}^n\omega_{ij}|f_i-f_j|-\mu^k(f\cdot v^k)\right\}$\\ 
2. {\bf Repeat}\\
\qquad i.$ f^{k+1}=g^{k+1}-\median(g^{k+1})$\\ 
\qquad ii. $f^{k+1}=\text{PseudoOrt}(g^{k+1},u_2)$\\
\quad{\bf Until} $\hat I( f^{k+1}) \ increases$\\
3. $v_i^{k+1}=\begin{cases}\sign(f_i^{k+1}), \qquad\qquad\text{if}\ f_i^{k+1}\neq 0,\\ -\dfrac{|f_+^{k+1}|-|f_-^{k+1}|}{|f_0^{k+1}|},\quad\text{if} \ f_i^{k+1}=0. \end{cases}$\\
4. $\mu^{k+1}=I(f^{k+1})$\\
{\bf Until} $\dfrac{|\mu^{k+1}-\mu^k|}{\mu^k}<\eps$\\ \hline
\end{tabular}
\end{center}

For the convergence of the Algorithm, we refer to \cite{BHb}, in particular Lemma 3.1, Theorems 3.1 and 4.1.

\begin{rem}
We focus on the $3$-clustering but these methods can be easily adapted to higher clustering.
As highlighted in Remark \ref{k_ort_rem}, for the $k$-clustering ($k>3$) we expect that \eqref{g_k_ort} are the conditions characterizing the $k$-th Cheeger constant. Therefore one can adapt the optimal tresholding of the second, the third, ... and the $k$-th eigenvector using the IPM as described before.

On the other hand, the proposed algorithm for $3$-clustering deeply rely on \cite{BHb}'s algorithms for $2$-clustering. So, a smart use of a combination of these algorithms could give very good approximation for the $k$-clustering with prescribed order $k>3$ .
\end{rem}

We modify the \cite{BHb}' code, to implement the described algorithms and methods on \texttrademark{MathLab} platform. The code is free downloadable at \url{https://github.com/GianpaoloPiscitelli/One_Spectral_3_Clustering}.

\section*{Acknowledgements}
This work has been partially supported by the MiUR-Dipartimenti di Eccellenza 2018-2022 grant \lq\lq Sistemi distribuiti intelligenti\rq\rq of Dipartimento di Ingegneria Elettrica e dell'Informazione \lq\lq M. Scarano\rq\rq, by the MiSE-FSC 2014-2020 grant \lq\lq SUMMa: Smart Urban Mobility
Management\rq\rq and by GNAMPA of INdAM.
We would also like to thank D.A. La Manna and V. Mottola for the helpful conversations during the  starting stage of this work.

\end{document}